\documentclass[11pt]{article}
\usepackage{amssymb,latexsym}
\usepackage{amsmath,amscd}
\usepackage{theorem}
\usepackage[all]{xy}
\title{\bf Poisson geometry and Morita equivalence}
\author{Henrique Bursztyn\thanks{{\tt henrique@math.toronto.edu},
    research partially supported by DAAD} \\[0.1cm]
        Department of Mathematics\\
    University of Toronto\\
        Toronto, Ontario M5S 3G3, Canada
        \\[0.2cm]
        Alan Weinstein \thanks{{\tt alanw@math.berkeley.edu}, research
    partially supported by NSF grant DMS-0204100
\newline \mbox{~~~~}MSC2000 Subject Classification Numbers: 53D17 (Primary), 58H05
16D90 (secondary)
\newline \mbox{~~~~}Keywords: Picard group, Morita equivalence, Poisson manifold,
symplectic groupoid, bimodule}
 \\[0.1cm]
         Department of Mathematics\\
        University of California, Berkeley\\
        CA, 94720-3840 USA\\and\\
Mathematical Sciences Research Institute\\17 Gauss Way\\Berkeley, CA
94720-5070 USA}
\date{}

\newcommand{\Lie}        {\mathcal L}

\newcommand{\id}         {{\mathrm {Id}}}

\newcommand{\Ker}        {{\mathrm {ker}}}

\newcommand{\Aut}        {{\mathrm {Aut}}}
\newcommand{\Inaut}      {{\mathrm {InnAut}}}
\newcommand{\Outaut}     {{\mathrm {OutAut}}}
\newcommand{\Ad}         {{\mathrm {Ad}}}
\newcommand{\gra}        {{\mathrm {graph}}}
\newcommand{\SP} [1]     {{\left\langle {{#1}} \right\rangle}}

\newcommand{\pr}         {{\mathrm{pr}}}
\newcommand{\Diff}       {\mathrm{Diff}}

\newcommand{\End}        {\mathrm{End}}
\newcommand{\st} [1]     {\scriptscriptstyle{{#1}}}

\newcommand{\inc}        {j}
\newcommand{\cent}       {\mathcal{Z}}
\newcommand{\proj}       {h}
\newcommand{\frakg}      {\mathfrak{g}}

\newcommand{\gstar}      {\mathfrak{g}^*}

\newcommand{\poismap}    {\psi}
\newcommand{\fib}        {\mathrm{Fib}}
\newcommand{\lcart}       {\lambda}
\newcommand{\rcart}       {\overline{\lambda}}

\newcommand{\grd}        {\mathcal{G}}

\newcommand{\Bi}[1]      {({#1})_{\frakg}}

\newcommand{\bim}         {\mathrm{Bim}}

\newcommand{\homom}       {\mathrm{Hom}}
\newcommand{\catalg}      {\mathsf{Alg}}
\newcommand{\catgr}       {\mathsf{LG}}

\newcommand{\SGrd}        {\mathsf{SG}}
\newcommand{\catpoiss}    {\mathsf{Poiss}}
\newcommand{\catCstar}    {\mathsf{C^*}}
\newcommand{\pt}          {\mathrm{pt}}

\newcommand{\A}          {\mathcal{A}}
\newcommand{\B}          {\mathcal{B}}
\newcommand{\C}          {\mathcal{C}}

\newcommand{\calo}       {\mathcal{O}}

\newcommand{\Cour}[1]      {[\![#1]\!]}

\newcommand{\Pic}         {\mathrm{Pic}}
\newcommand{\SPic}        {\mathrm{SPic}}

\newcommand{\X}          {\mathrm{X}}
\newcommand{\Y}          {\mathrm{Y}}

\newcommand{\AXB}        {{_{\scriptscriptstyle{\A}}{\X}_{\scriptscriptstyle{\B}}}}

\newcommand{\BYA}        {{_{\scriptscriptstyle{\B}}{\Y}_{\scriptscriptstyle{\A}}}}
\newcommand{\BYC}        {{_{\scriptscriptstyle{\B}}{\Y}_{\scriptscriptstyle{\C}}}}

\newcommand{\GXH}        {{_{\scriptscriptstyle{G}}{\X}}_{\scriptscriptstyle{H}}}
\newcommand{\HYK}        {{_{\scriptscriptstyle{H}}{\Y}}_{\scriptscriptstyle{K}}}
\newcommand{\PSP}        {P_1\stackrel{J_1}{\leftarrow}S\stackrel{J_2}{\rightarrow}P_2}
\newcommand{\PSbarP}     {P_2\stackrel{J_2}{\leftarrow}\overline{S}\stackrel{J_1}{\rightarrow}P_1}

\newcommand{\PGPst}      {P\stackrel{t}{\leftarrow}\grd\stackrel{s}{\rightarrow}P}
\newcommand{\PSPprime}   {P_2\stackrel{J_2'}{\leftarrow}S'\stackrel{J_3'}{\rightarrow}P_3}

\newcommand{\reals}      {{\mathbb R}}
\newcommand{\integers}   {{\mathbb Z}}

\newcommand{\repC}      {{\mathrm{Herm}}}

\newtheorem{lemma} {Lemma} [section]
\newtheorem{proposition} [lemma] {Proposition}

\newtheorem{theorem} [lemma] {Theorem}

\theorembodyfont{\rm} 

\newtheorem{example}[lemma] {Example}

\newtheorem{remark}[lemma]{Remark}

\newenvironment{proof}{{\sc Proof:}}{{\hspace*{\fill} $\square$\\}}

\newenvironment{exercise}{\medskip \begin{center}
\begin{minipage}{5.5in} \footnotesize \noindent {\bf Exercise \\ }}{\end{minipage}
\end{center}
\medskip}

\numberwithin{equation}{section}

\begin{document}

\maketitle



\tableofcontents

\section{Introduction}

Poisson geometry is a ``transitional'' subject between noncommutative
algebra and differential geometry (which could be seen as the study of
a very special class of commutative algebras).  The physical
counterpart to this transition is the
correspondence principle linking quantum to classical mechanics.

The main purpose of these notes is to present an aspect of Poisson
geometry which is inherited from the noncommutative side: the notion
of Morita equivalence, including the ``self-equivalences'' known as
Picard groups.

In algebra, the importance of Morita equivalence lies in the fact that
Morita equivalent algebras have, by definition, equivalent categories
of modules.  From this it follows that many other invariants, such as
cohomology and deformation theory, are shared by all Morita equivalent
algebras.  In addition, one can sometimes understand the
representation theory of a given algebra by analyzing that of a
simpler representative of its Morita equivalence class.  In Poisson
geometry, the role of ``modules'' is played by Poisson maps from
symplectic manifolds to a given Poisson  manifold.  The
simplest such maps are the inclusions of symplectic leaves, and indeed
the structure of the leaf space is a Morita invariant.  (We will see
that this leaf space sometimes has a more rigid structure than one
might expect.)

The main theorem of algebraic Morita theory is that Morita
equivalences are implemented by bimodules.  The same thing turns out
to be true in Poisson geometry, with the proper geometric definition
of ``bimodule''.

Here is a brief outline of what follows this introduction.

Section 2 is an introduction to Poisson geometry and some of its
recent generalizations, including Dirac geometry and ``twisted''
Poisson geometry in the presence of a ``background'' closed 3-form.
Both of these generalizations are used simultaneously to get a geometric
understanding of new notions of symmetry of growing importance in
mathematical physics, especially with background 3-forms arising
throughout string theory (in the guise of the more familiar closed 2- forms
on spaces of curves).

In Section 3, we review various flavors of the algebraic theory of
Morita equivalence in a way which transfers easily to the geometric
case.  In fact, some of our examples come from geometry: algebras of
smooth functions.  Others come from the quantum side: operator
algebras.

Section 4 is the heart of these notes, a presentation of the geometric
Morita theory of Poisson manifolds and the closely related Morita
theory of symplectic groupoids.  We arrive at this theory via
the Morita theory of Lie groupoids in general.

In Section 5, we attempt to remedy a defect in the theory of Chapter
4.  Poisson manifolds with equivalent (even isomorphic) representation
categories may not be Morita equivalent.  We introduce refined versions
of the representation category (some of which are not really
categories!) which do determine the Morita equivalence class.  Much of
the material in this section is new and has not yet appeared in
print.  (Some of it is based on discussions which came after the PQR
Euroschool where this course was presented.)

Along the way, we comment on a pervasive problem in the geometric
theory.  Many constructions involve forming the leaf space of a
foliation, but these leaf spaces are not always manifolds.  We make
some remarks about the use of differentiable stacks as a language for
admitting pathological leaf spaces into the world of smooth geometry.

\noindent {\bf Acknowledgements:}

We would like to thank all the organizers and participants at the
Euroschool on Poisson Geometry, Deformation Quantization, and
Representation Theory for the opportunity to present this short
course, and for their feedback at the time of the School. We also thank
Stefan Waldmann for his comments on the manuscript.

H.B. thanks Freiburg University for its hospitality while part of this work was being done.

\section{Poisson geometry and some  generalizations}

\subsection{Poisson manifolds}
Let $P$ be a smooth manifold. A \textbf{Poisson structure} on $P$
is an $\reals$-bilinear Lie bracket $\{\cdot,\cdot\}$ on $C^\infty(P) $ satisfying the Leibniz rule
\begin{equation}\label{eq:leibniz}
\{f,gh\} = \{f,g\}h + g\{f,h\}, \;\; \mbox{ for all } f,g,h \in
C^\infty(P).
\end{equation}
A \textbf{Poisson algebra} is an associative algebra which is
also a Lie algebra so that the associative multiplication and the
Lie bracket are related by \eqref{eq:leibniz}.

For a function $f \in C^\infty(P)$, the derivation $X_f =
\{f,\cdot \}$ is called the \textbf{hamiltonian vector field} of
$f$. If $X_f=0$, we call $f$ a \textbf{Casimir function} (see
Remark \ref{rmk:casimir}). It follows from \eqref{eq:leibniz} that
there exists a bivector field $\Pi \in \mathcal{X}^2(P) =
\Gamma(\bigwedge^2 TP)$ such that
$$
\{f,g\}= \Pi(df,dg);
$$
the Jacobi identity for $\{\cdot,\cdot\}$ is equivalent to the condition
$[\Pi,\Pi]=0$, where $[\cdot,\cdot]$ is the Schouten- Nijenhuis
bracket, see e.g.  \cite{Vais}.

In local coordinates $(x_1,\cdots, x_n)$, the tensor $\Pi$ is determined by the matrix
\begin{equation}\label{eq:matrix}
\Pi_{ij} (x) = \{x_i,x_j\}.
\end{equation}
If this matrix is invertible at each $x$, then $\Pi$ is
called nondegenerate or \textbf{symplectic}. In this case, the
local matrices $(\omega_{ij}) = (-\Pi_{ij})^{-1}$ define a global 2-form $\omega
\in \Omega^2(P)=\Gamma(\bigwedge^2T^*P)$, and the condition $[\Pi,\Pi]=0$ is
equivalent to $d\omega=0$.

\begin{example}(\textit{Constant Poisson structures})\label{ex:constant}

Let $P=\mathbb{R}^n$, and suppose that the $\Pi_{ij} (x)$ are constant.
By a linear change of coordinates, one can find new coordinates
$(q_1,\cdots,q_k,p_1,\cdots,p_k,e_1,\cdots,e_l)$, $2k+l=n$, so that
$$
\Pi=\sum_i \frac{\partial}{\partial q_i}\wedge
\frac{\partial}{\partial p_i}.
$$
In terms of the bracket, we have
$$
\{f,g\}=\sum_i\left(\frac{\partial f}{\partial q_i}\frac{\partial g}{\partial p_i}
- \frac{\partial f}{\partial p_i}\frac{\partial g}{\partial q_i}\right)
$$
which is the original Poisson bracket in mechanics. In this example,
all the coordinates $e_i$ are Casimirs.
\end{example}

\begin{example} (\textit{Poisson structures on $\mathbb{R}^2$})\label{ex:poissonR2}

Any smooth function $f: \mathbb{R}^2 \to \mathbb{R}$ defines a
Poisson structure in $\mathbb{R}^2=\{(x_1,x_2)\}$ by
$$
\{x_1,x_2\} := f(x_1,x_2),
$$
and every Poisson structure on $\reals^2$ has this form.
\end{example}

\begin{example}\label{ex:liepoisson} (\textit{Lie-Poisson structures})

An important class of Poisson structures are the linear
ones.  If $P$ is a (finite-dimensional)
vector space $V$ considered as a manifold,  with linear
coordinates $(x_1,\cdots,x_n)$, a linear Poisson structure is determined
by constants $c_{ij}^k$ satisfying
\begin{equation}\label{eq:linear}
\{x_i,x_j\}=\sum_{k = 1}^n c_{ij}^k x_k.
\end{equation}
(We may assume that $c_{ij}^k=-c_{ji}^k.$)
Such Poisson structures are usually called \textbf{Lie-Poisson structures},
since the Jacobi identity for the Poisson bracket implies that the
$c_{ij}^k$ are the structure constants of a Lie algebra
$\mathfrak{g}$, which may be identified in a natural way with $V^*$.
(Also, these Poisson structures were originally introduced by Lie \cite{li:theorie} himself.)
Note that we may also identify $V$ with $ \mathfrak{g}^*$.
Conversely, any Lie
algebra $\mathfrak{g}$ with structure constants $c_{ij}^k$ defines
by \eqref{eq:linear} a linear Poisson structure on $\mathfrak{g}^*$.
\end{example}

\begin{remark}(\textit{Casimir functions})\label{rmk:casimir}

Deformation quantization of the Lie-Poisson structure on $\frakg^*$,
see e.g. \cite{be:remarks,gu:explicit},
leads to the universal enveloping algebra $U(\frakg)$.  Elements of
the center of $U(\frakg)$ are known as Casimir elements (or Casimir
operators, when a representation of $\frakg$ is extended to a
representation of $U(\frakg)$).  These correspond to the center of the
Poisson algebra of functions on $\frakg^*$, hence, by extension, the
designation ``Casimir functions'' for the center of any Poisson algebra.
\end{remark}

\subsection{Dirac structures}\label{subsec:dirac}

We now introduce a simultaneous generalization of Poisson structures and closed 2-forms. (We will often
refer to closed 2-forms as \textbf{presymplectic}.)

Each 2-form $\omega$ on $P$ corresponds to a bundle
map
\begin{equation}\label{eq:bundle1}
\widetilde{\omega}:TP \to T^*P, \;\;\;
\widetilde{\omega}(v)(u)=\omega(v,u).
\end{equation}
Similarly, for a bivector field $\Pi \in \mathcal{X}^2(P)$, we
define the bundle map
\begin{equation}\label{eq:bundle2}
\widetilde{\Pi}:T^*P \to TP, \;\;\;
\beta(\widetilde{\Pi}(\alpha))=\Pi(\alpha,\beta).
\end{equation}
The matrix representing $\widetilde{\Pi}$ in the bases $(dx_i)$
and $(\partial/\partial x_i)$ corresponding to local coordinates
induced by coordinates $(x_1,\ldots,x_n)$ on $P$ is, up to a sign,
just \eqref{eq:matrix}. So bivector fields (or 2-forms) are
nondegenerate if and only if the associated bundle maps are
invertible.

By using the maps in \eqref{eq:bundle1} and \eqref{eq:bundle2},
we can describe both closed 2-forms and Poisson bivector fields as subbundles of
$TP\oplus T^*P$: we simply consider the graphs
$$
L_{\omega}:= \gra(\widetilde{\omega}),\;\; \mbox{ and }\;\;
L_{\Pi}:= \gra(\widetilde{\Pi}).
$$
To see which subbundles of $TP\oplus T^*P$ are of this form, we introduce the
following canonical structure on $TP\oplus T^*P$:

\begin{itemize}
\item[1)] The symmetric bilinear form $\SP{\cdot,\cdot}_+:TP \oplus T^*P \to \mathbb{R}$,
\begin{equation}\label{eq:symbr}
\SP{(X,\alpha), (Y,\beta)}_+:=\alpha(Y) + \beta(X).
\end{equation}
\item[2)] The bracket $\Cour{\cdot,\cdot}:\Gamma(TP\oplus T^*P)\times
\Gamma(TP\oplus T^*P) \to \Gamma(TP\oplus T^*P)$,
\begin{equation} \label{eq:coubr}
\Cour{(X,\alpha), (Y,\beta)}:= ([X,Y], \mathcal{L}_X\beta -i_Yd\alpha).
\end{equation}
\end{itemize}

\begin{remark}(\textit{Courant bracket})

The bracket \eqref{eq:coubr} is the non-skew-symmetric version,
introduced in \cite{LWX} (see also \cite{SeWe01}),  of T.~Courant's
original bracket \cite{Cou90}. The bundle $TP\oplus T^*P$ together
with the brackets \eqref{eq:symbr} and \eqref{eq:coubr} is an
example of a \textbf{Courant algebroid} \cite{LWX}.
\end{remark}

Using the brackets \eqref{eq:symbr} and \eqref{eq:coubr}, we have
the following result \cite{Cou90}:

\begin{proposition}\label{prop:dirac}
A subbundle $L\subset TP\oplus T^*P$ is of the form
$L_{\Pi}=\gra(\widetilde{\Pi})$ (resp. $L_{\omega} =
\gra(\widetilde{\omega})$) for a bivector field $\Pi$ (resp.
2-form $\omega$)
 if and only if
\begin{itemize}

\item[i)] $TP\cap L=\{0\}$ (resp. $L\cap T^*P=\{0\}$) at all
points of $P$;

\item[ii)] $L$ is maximal
  isotropic with respect to
$\SP{\cdot,\cdot}_+$;
\end{itemize}
furthermore, $[\Pi,\Pi]=0$ (resp. $d\omega=0$) if and only if
\begin{itemize}
\item[iii)] $\Gamma(L)$ is closed under the Courant bracket
\eqref{eq:coubr}.

\end{itemize}
\end{proposition}

Recall that $L$ being isotropic with respect to $\SP{\cdot,\cdot}_+$ means that,
at each point of $P$,
$$
\SP{(X,\alpha),(Y,\beta)}_+=0
$$
whenever $(X,\alpha), (Y,\beta) \in L$.  Maximality is equivalent
to the dimension condition $\mathrm{rank}(L)= \mathrm{dim}(P)$.

A \textbf{Dirac structure} on $P$ is a subbundle $L \subset TP
\oplus T^*P$ which is maximal isotropic with respect to $\SP{\cdot,\cdot}_+$
and whose sections are closed under the Courant bracket \eqref{eq:coubr};
in other words, a Dirac structure satisfies conditions $ii)$ and $iii)$ of
Prop.~\ref{prop:dirac} but is not necessarily the graph associated to a
bivector field or 2-form.

If $L$ satisfies only $ii)$, it is called an \textbf{almost Dirac
structure}, and we refer to $iii)$ as the \textbf{integrability
condition} of a Dirac structure. The next example illustrates
these notions in another situation.

\begin{example}(\textit{Regular foliations})

Let $F \subseteq TP$ be a subbundle, and let $F^\circ \subset
T^*P$ be its annihilator. Then $L=F\oplus F^\circ$ is an almost
Dirac structure; it is a Dirac structure if and only if $F$
satisfies the Frobenius condition
$$
[\Gamma(F),\Gamma(F)] \subset \Gamma(F).
$$
So regular foliations are examples of Dirac structures.
\end{example}

\begin{example}\label{ex:lineardirac}({\textit{Vector Dirac structures}})

If $V$ is a finite-dimensional real vector space, then a
\textbf{vector Dirac structure} on $V$ is a subspace $L \subset
V\oplus V^*$ which is maximal isotropic with respect to the
symmetric pairing \eqref{eq:symbr}.\footnote{Vector Dirac structures
  are sometimes called ``linear Dirac structures," but we will eschew
  this name  to avoid confusion with linear (i.e. Lie-) Poisson
  structures.  (See Example \ref{ex:liepoisson})}

Let $L$ be a vector Dirac structure on $V$. Let
$\pr_1:V \oplus V^* \to V$ and $\pr_1:V\oplus V^* \to V^*$ be the
canonical projections, and consider the subspace
$$
R := \pr_1(L) \subseteq V.
$$
Then $L$ induces a skew-symmetric bilinear form $\theta$ on $R$
defined by
\begin{equation}\label{eq:diracform}
\theta(X,Y):=\alpha(Y),
\end{equation}
where $X,Y \in R$ and $\alpha \in V^*$ is such that $(X,\alpha)
\in L$.

\begin{exercise}
Show that $\theta$ is well defined, i.e., \eqref{eq:diracform} is
independent of the choice of $\alpha$.
\end{exercise}

Conversely, any pair $(R,\theta)$, where $R\subseteq V$ is a
subspace and $\theta$ is a skew-symmetric bilinear form on $R$,
defines a vector Dirac structure by
\begin{equation}\label{eq:Linduced}
L:=\{(X,\alpha),\;\; X\in R,\; \alpha \in V^* \mbox{ with }
\alpha|_R = i_{X}\theta\}.
\end{equation}

\begin{exercise}
Check that $L$ defined in \eqref{eq:Linduced} is a vector Dirac
structure on $V$ with associated subspace $R$ and bilinear form
$\theta$.
\end{exercise}
\end{example}

Example \ref{ex:lineardirac} indicates a simple way in which vector Dirac structures can be
restricted to subspaces.

\begin{example}\label{ex:linearrest} (\textit{Restriction of Dirac structures to subspaces})

Let $L$ be a vector Dirac structure on  $V$, let
$W \subseteq V$ be a subspace, and consider the pair $(R,\theta)$
associated with $L$. Then $W$
inherits the vector Dirac structure $L_W$ from $L$
defined by the pair
$$
R_W:= R\cap W, \; \mbox{ and } \; \theta_W:=\iota^*\theta,
$$
where $\iota:W \hookrightarrow V$ is the inclusion map.

\begin{exercise}
Show that there is a canonical isomorphism
\begin{equation}\label{eq:canisom}
L_W \cong \frac{L\cap(W\oplus V^*)}{L\cap W^\circ}.
\end{equation}
\end{exercise}

\end{example}

Let $(P,L)$ be a Dirac manifold, and let $\iota:N \hookrightarrow P$ be  a submanifold.
The construction in Example \ref{ex:linearrest}, when applied to
$T_x N\subseteq T_x P$ for all $x\in P$,
defines a maximal isotropic ``subbundle'' $L_N \subset TN \oplus T^*N$. The problem
is that $L_N$ may not be a continuous family of subspaces.  When  $L_N$ {\em is} a continuous family, it is a smooth bundle which then automatically
satisfies the integrability condition \cite[Cor.~3.1.4]{Cou90},
so $L_N$ defines a Dirac structure on $N$.

The next example is a special case of this construction and is one of the original motivations
for the study of Dirac structures; it illustrates the connection
between Dirac structures and ``constraint submanifolds'' in classical mechanics.

\begin{example}(\textit{Momentum level sets})\label{ex:momlevel}

Let $J:P\to \mathfrak{g}^*$ be the momentum \footnote{The term
  ``moment'' is frequently used instead of ``momentum'' in this
  context.  In this paper, we will follow the convention, introduced
  in \cite{MiWe}, that ``moment'' is used only in connection with
  groupoid actions.  As we will see (e.g. in Example \ref{ex:hamilt}),
many momentum maps, even for
  ``exotic'' theories, are moment maps as well.}
map for a hamiltonian action of a Lie
group $G$ on a Poisson manifold $P$ \cite{MaRa}.
Let $\mu \in \mathfrak{g}^*$ be a regular value for
$J$, let $G_{\mu}$ be the isotropy group at $\mu$ with respect to
the coadjoint action, and consider
$$
Q=J^{-1}(\mu) \hookrightarrow P.
$$
At each point $x\in Q$, we have a vector Dirac structure on $T_xQ$ given by
\begin{equation}\label{eq:levelsetdirac}
(L_Q)_x := \frac{L_x\cap(T_xQ\oplus T^*_xP)}{L_x\cap T_xQ^\circ}.
\end{equation}
To show  that $L_Q$ defines a smooth bundle, it suffices to verify
that $L_x\cap T_xQ^\circ$ has constant dimension. (Indeed, if this is the case,
then $L_x\cap(T_xQ\oplus T^*_xP)$ has constant dimension as well, since the quotient
$L_x\cap(T_xQ\oplus T^*_xP)/{L_x\cap T_xQ^\circ}$ has constant dimension, and this
insures that all bundles are smooth.)
A direct computation shows that $L_x\cap T_xQ^\circ$ has constant dimension if and only if
the stabilizer groups of the $G_\mu$-action on $Q$ have constant dimension, which happens
whenever the $G_\mu$-orbits on $Q$ have constant dimension
(for instance, when the action of $G_\mu$ on $Q$ is locally free).
In this case, $L_Q$ is a Dirac structure on $Q$.
\end{example}

We will revisit this example in Section \ref{subsec:diracmaps}.

\begin{remark}(\textit{Complex Dirac structures and generalized complex geometry})\label{rem:complex}

Using the natural extensions of the symmetric form
\eqref{eq:symbr} and the Courant bracket \eqref{eq:coubr} to
$(TP\oplus T^*P)\otimes \mathbb{C}$, one can define a
\textbf{complex Dirac structure} on a manifold $P$ to be a maximal
isotropic \textit{complex} subbundle $L \subset (TP\oplus
T^*P)\otimes \mathbb{C}$ whose sections are closed under the
Courant bracket. If a complex Dirac structure $L$ satisfies the condition
\begin{equation}\label{eq:complexcond}
L\cap \overline{L}=\{0\}
\end{equation}
at all points of $P$ (here $\overline{L}$ is the complex conjugate
of $L$), then it is called a \textbf{generalized complex
structure}; such structures were introduced in
\cite{Gualt,Hitchin} as a common generalization of complex and
symplectic structures.

To see how complex structures fit into this picture,
note that an almost complex structure
$J:TP\to TP$ defines a maximal isotropic subbundle
$L_J \subset (TP\oplus T^*P)\otimes \mathbb{C}$ as the
$i$-eigenbundle of the map
$$
(TP\oplus T^*P)\otimes \mathbb{C} \to (TP\oplus T^*P)\otimes \mathbb{C},
\;\;\; (X,\alpha) \mapsto
(-J(X),J^*(\alpha)).
$$
The bundle $L_J$ completely characterizes $J$, and
it satisfies \eqref{eq:complexcond}; moreover, $L_J$
satisfies the integrability condition of a Dirac structure
if and only if $J$ is a complex structure.

Similarly, a symplectic structure $\omega$ on $P$ can be seen
as a generalized complex structure through the bundle
$L_{\omega,{\mathbb{C}}}$,
defined as the $i$-eigenbundle of the map
$$
(TP\oplus T^*P)\otimes \mathbb{C} \to (TP\oplus T^*P)\otimes \mathbb{C},
\;\;\;
(X,\alpha) \mapsto (\widetilde{\omega}(X),-\widetilde{\omega}^{-1}(\alpha)).
$$

Note that, by \eqref{eq:complexcond}, a generalized complex structure is never
the complexification of a real Dirac structure.  In particular, for a
symplectic structure $\omega$, $L_{\omega,\mathbb{C}}$ is \textit{not} the
  complexification of the real Dirac structure $L_\omega$ of
  Proposition \ref{prop:dirac}.
\end{remark}

\subsection{Twisted structures}\label{subsec:twisted}

A ``background" \textit{closed} 3-form $\phi \in \Omega^3(P)$ can be used to
``twist'' the geometry of $P$ \cite{KlStr,Park}, leading to a modified notion of
 Dirac structures \cite{SeWe01}, and in particular of Poisson structure. The
key point is to use $\phi$ to alter the ordinary Courant bracket
\eqref{eq:coubr} as follows:
\begin{equation}\label{eq:twistcour}
\Cour{(X,\alpha),(Y,\beta)}_{\phi}:=([X,Y], \mathcal{L}_X\beta -
i_Yd\alpha + \phi(X,Y,\cdot)).
\end{equation}
We now simply repeat the definitions in Section \ref{subsec:dirac}
replacing \eqref{eq:coubr} by the \textbf{$\phi$-twisted Courant
bracket} \eqref{eq:twistcour}.

A \textbf{$\phi$-twisted Dirac structure} on $P$ is a subbundle $L
\subset TP\oplus T^*P$ which is maximal isotropic with respect to
$\SP{\cdot,\cdot}_+$ \eqref{eq:symbr} and for which
\begin{equation}
\Cour{\Gamma(L),\Gamma(L)}_{\phi} \subseteq \Gamma(L).
\end{equation}
With this new integrability condition, one can check that the graph of a bivector field $\Pi$ is a
$\phi$-twisted Dirac structure if and only if
$$
\frac{1}{2}[\Pi,\Pi]=\wedge^3\widetilde{\Pi}(\phi);
$$
such bivector fields are called \textbf{$\phi$-twisted Poisson
structures}. Similarly,  the graph of a 2-form $\omega$ is a $\phi$-twisted Dirac structure
if and only if
$$
d\omega + \phi =0,
$$
in which case $\omega$ is called a \textbf{$\phi$-twisted
presymplectic structure}.

\begin{remark}(\textit{Terminology})

The term ``twisted Dirac structure'' and its cousins represent a
certain abuse of terminology, since it is not the Dirac (or Poisson,
etc.) structure which is twisted, but rather the \textit{notion} of
Dirac structure.  Nevertheless, we have chosen to stick to this
terminology, rather than the alternative ``Dirac structure with
background'' \cite{K-S03}, because it is consistent with such existing
terms as ``twisted sheaf'', and because the alternative terms lead to
some awkward constructions.
\end{remark}

\begin{example}\label{ex:cartan-dirac}(\textit{Cartan-Dirac structures on Lie groups})

Let $G$ be a Lie group whose Lie algebra $\mathfrak{g}$ is equipped with a
nondegenerate adjoint-invariant symmetric bilinear form
$\Bi{\cdot,\cdot}$, which we use to identify $TG$ and
$T^*G$. In $TG\oplus TG\sim TG\oplus T^*G$, we consider the maximal isotropic
subbundle
\begin{equation}\label{eq:cartandirac}
L_G := \{ (v_r-v_l,\frac{1}{2}(v_r+v_l)),\;\; v \in \mathfrak{g}\},
\end{equation}
where $v_r$ and $v_l$ are the right and left invariant vector
fields corresponding to $v$. One can show that $L_G$ is a
$\phi^G$-twisted Dirac structure, where $\phi^G$ is the bi-invariant
Cartan 3-form on $G$, defined on Lie algebra elements by
$$
\phi^G(u,v,w)=\frac{1}{2}\Bi{u,[v,w]}.
$$
We call $L_G$ the \textbf{Cartan-Dirac structure} on $G$ associated
with $\Bi{\cdot,\cdot}$. Note that $L_G$ is of the form $L_{\Pi}$
only at points $g$ for which $\Ad_g+1$ is invertible, see also Example \ref{ex:folCD}.

These Dirac structures are closely related to the theory of quasi-hamiltonian
spaces and group-valued momentum maps \cite{AMM,BCWZ,Xu03}, as well as to
quasi-Poisson manifolds \cite{AKM,BuCr}.
\end{example}

\subsection{Symplectic leaves and local structure of Poisson manifolds}\label{subsec:local}

If $\Pi$ is a symplectic Poisson structure on $P$, then
Darboux's theorem asserts that, around each point of $P$, one can
find coordinates $(q_1,\cdots,q_k,p_1,\cdots,p_k)$ such that
$$
\Pi=\sum_i \frac{\partial}{\partial q_i}\wedge
\frac{\partial}{\partial p_i}.
$$
The corresponding symplectic form $\omega$ is
$$
\omega=\sum_i dq_i \wedge dp_i.
$$

In general, the image of the bundle map \eqref{eq:bundle2},
$\widetilde{\Pi}(T^*P) \subseteq TP$, defines an integrable
singular distribution on $P$; in other words, $P$ is a disjoint
union of ``leaves'' $\mathcal{O}$ satisfying $T_x\mathcal{O} =
\widetilde{\Pi}(T_x^*P)$ for all $x \in P$. The leaf $\mathcal{O}$
through $x$ can be described as the points which can be reached
from $x$ through piecewise hamiltonian paths.

If $\widetilde{\Pi}$ has locally constant rank, we call the
Poisson structure $\Pi$ \textbf{regular}, in which case it
defines a foliation of $P$ in the ordinary sense. Note that this
is always the case on an open dense subset of $P$, called the
\textbf{regular part}.

The local structure of a Poisson manifold $(P,\Pi)$ around a
regular point is given the Lie-Darboux theorem: If $\Pi$ has
constant rank $k$ around a given point, then there exist
coordinates $(q_1,\ldots,q_k,p_1,\ldots,p_k,e_1,\ldots,e_l)$ such
that
$$
\{q_i,p_j\}=\delta_{ij},\; \mbox{ and }\; \{q_i,q_j\}=\{p_i,p_j\}=
\{q_i,e_j\}=\{p_i,e_j\}=0.
$$
Thus, the local structure of a regular Poisson manifold is determined by that of the vector Poisson structures on any of its tangent spaces (in a given connected component).

In the general case, we have the local splitting theorem
\cite{We83}:
\begin{theorem}\label{thm:split}
Around any point $x_0$ in a Poisson manifold $P$, there exist
coordinates
$$(q_1,\ldots,q_k,p_1,\ldots,p_k,e_1,\ldots,e_l), \;\;\; (q,p,e)(x_0)=(0,0,0),$$
such that
$$
\Pi= \sum_{i=1}^k\frac{\partial }{\partial
q_i}\wedge\frac{\partial}{\partial p_i} +
\frac{1}{2}\sum_{i,j=1}^l \eta_{ij}(e)\frac{\partial}{\partial
e_i}\wedge\frac{\partial}{\partial e_j}
$$
and $\eta_{ij}(0) = 0$.
\end{theorem}

The splitting of Theorem
\ref{thm:split} has a symplectic factor associated with the
coordinates $(q_i,p_i)$
and a totally degenerate factor (i.e., with all Poisson brackets
vanishing at $e=0$) associated with the coordinates $(e_j)$.
The symplectic factor may be identified with an open subset of the leaf $\mathcal{O}$
through $x_0$; patching them together defines
a symplectic structure on each leaf of the foliation determined by $\Pi$.
So $\Pi$ canonically defines a singular \textbf{symplectic foliation} of
$P$.
The totally degenerate factor in the local splitting is
well-defined up to isomorphism. Its isomorphism class is the same at all points in a given symplectic leaf, so one refers to the
totally degenerate factor as the \textbf{transverse structure} to
$\Pi$ along the leaf.

\begin{example} (\textit{Symplectic leaves of Poisson structures on $\mathbb{R}^2$})

Let $f: \mathbb{R}^2 \to \mathbb{R}$ be a smooth function, and
let us consider the Poisson structure on  $\mathbb{R}^2=\{(x_1,x_2)\}$ defined by
$$
\{x_1,x_2\} := f(x_1,x_2).
$$
The connected components of the set where
 $f(x_1,x_2)\neq 0$ are the
2-dimensional symplectic leaves; in the set where $f$ vanishes,
each point is a symplectic leaf.
\end{example}

\begin{example}\label{ex:coadleaves} (\textit{Symplectic leaves of Lie-Poisson structures})

Let us consider $\mathfrak{g}^*$, the dual of the Lie algebra $\mathfrak{g}$,
equipped with its Lie-Poisson structure, see Example \ref{ex:liepoisson}.
The symplectic leaves are just the coadjoint orbits for any connected group with
Lie algebra $\mathfrak{g}$.  Since $\{0\}$ is always an orbit,
a Lie-Poisson structure
is not regular unless $\mathfrak{g}$ is abelian.
\end{example}

\begin{exercise}
Describe the symplectic leaves in the duals of $\mathfrak{su}(2)$,
$\mathfrak{sl}(2,\reals)$ and $\mathfrak{a}(1)$ (nonabelian
2-dimensional Lie algebra).
\end{exercise}

\begin{remark}(\textit{Linearization problem})

By linearizing at $x_0$ the functions $\eta_{ij}$  in Theorem
\ref{thm:split}, we can write
\begin{equation}\label{eq:linearized}
\{e_i,e_j\}=\sum_{k}c_{ij}^k e_k + O(e^2),
\end{equation}
and it turns out that $c_{ij}^k$ define a Lie-Poisson structure on
the normal space to the symplectic leaf at $x_0$. The
\textbf{linearization problem} consists of determining whether one
can choose suitable ``transverse'' coordinates $(e_1,\ldots,e_l)$
with respect to which $O(e^2)$ in \eqref{eq:linearized} vanishes.
For example, if the Lie algebra structure on the conormal bundle
to a symplectic leaf determined by the linearization of $\Pi$ at a
point $x_0$  is semi-simple and of compact type, then $\Pi$ is
linearizable around $x_0$ through a smooth change of coordinates.
The first proof of this theorem, due to Conn \cite{Conn}, used
many estimates and a ``hard'' implicit function theorem.  A
``soft'' proof, using only the sort of averaging usually
associated with compact group actions (but for groupoids instead
of groups), has recently been announced by Crainic and Fernandes \cite{CrFe3}.
There is also a ``semilocal'' problem of linearization in the
neighborhood of an entire symplectic leaf.   This problem was first
addressed by Vorobjev \cite{Vo}, with further developments by Davis
and Wade \cite{DaWa}.
For overviews of linearization and more general normal form questions,
we refer to the article of Fernandes and Monnier \cite{FeMo} and the
forthcoming book of Dufour and Zung \cite{DuZu}.
\end{remark}


\subsection{Presymplectic leaves and Poisson quotients of Dirac manifolds}
\label{subsec:quotients}

Let $\pr_1: TP\oplus T^*P \to TP$ and $\pr_2:TP\oplus T^*P \to
T^*P$ be the canonical projections. If $L \subset TP \oplus T^*P$
is a (twisted) Dirac structure on $P$, then
\begin{equation}\label{eq:diracdist}
\pr_1(L) \subseteq TP
\end{equation}
defines a singular distribution on $P$. Note that if $L=L_{\Pi}$ for a
Poisson structure $\Pi$, then $\pr_1(L)=\widetilde{\Pi}(T^*P)$, so
this distribution coincides with the one defined by $\Pi$, see Section \ref{subsec:local}.
It turns out that the integrability condition for (twisted) Dirac structures
guarantees that \eqref{eq:diracdist} is integrable in general, so
a (twisted) Dirac structure $L$ on $P$ determines a decomposition of $P$ into
leaves $\mathcal{O}$ satisfying
$$
T_x\mathcal{O} = \pr_1(L)_x
$$
at all $x\in P$.

Just as leaves of foliations associated with Poisson structures
carry symplectic forms, each leaf of a (twisted) Dirac manifold $P$ is
naturally equipped with a (twisted) presymplectic 2-form $\theta$:
at each $x \in P$,
$\theta_x$ is the bilinear form defined in \eqref{eq:diracform}. These
forms fit together into a smooth leafwise 2-form, which is
nondegenerate on the leaves just when $L$ is a (twisted) Poisson
structure. If $L$ is twisted by $\phi$, then $\theta$ is twisted by
the pull back of $\phi$ to each leaf.

\begin{remark}\label{rem:algebroids}(\textit{Lie algebroids})

The fact that $\pr_1(L) \subseteq TP$ is an integrable singular
distribution is a consequence of a more general fact: the
restriction of the Courant bracket $\Cour{\cdot,\cdot}_{\phi}$ to
$\Gamma(L)$ defines a Lie algebra bracket making $L\to P$ into a
\textit{Lie algebroid} with anchor $\pr_1|_L$, and the image of the
anchor of any Lie algebroid is always an integrable distribution
(its leaves are also called \textbf{orbits}). We refer to \cite{SilWein99,MM03}
for more on Lie algebroids.
\end{remark}

\begin{example}(\textit{Presymplectic leaves of Cartan-Dirac structures})\label{ex:folCD}

Let $L_G$ be a  Cartan-Dirac structure on $G$ with respect to $\Bi{\cdot,\cdot}$, see
\eqref{eq:cartandirac}. Then the associated distribution on $G$ is
$$
\pr_1(L_G)=\{v_r-v_l,\; v \in \mathfrak{g}\}.
$$
Since vector fields of the form $v_r-v_l$ are infinitesimal
generators of the action of $G$ on itself by conjugation, it
follows that the twisted presymplectic leaves of $L_G$ are the connected
components of the conjugacy classes of $G$.
With $v_G=v_r-v_l$, the corresponding twisted presymplectic forms can
be written as
\begin{equation}\label{eq:GHJW}
\theta_g(v_G,w_G):=\frac{1}{2}\Bi{(\Ad_{g^{-1}}-\Ad_{g})v,w},
\end{equation}
at $g \in G$. These 2-forms were introduced in \cite{GHJW} in the
study of the symplectic structure of certain moduli spaces.
They are analogous to the
Kostant-Kirillov-Souriau symplectic forms on coadjoint orbits,
although they are neither nondegenerate nor closed: $\theta_g$ is
degenerate whenever $1 + \Ad_g$ is not invertible, and, on a
conjugacy class $\iota:\mathcal{O}\hookrightarrow G$, $d\theta =
-\iota^*\phi^G$.

Just as the symplectic forms along coadjoint orbits
on the dual of a Lie algebra are associated with Lie-Poisson
structures, the 2-forms \eqref{eq:GHJW} along conjugacy classes of
a Lie group are associated with Cartan-Dirac structures.
\end{example}

For any $\phi$-twisted Dirac structure $L$, the (topologically)
closed family of subspaces
$TP\cap L=\ker(\theta)$ in $TP$ is called the \textbf{characteristic
distribution} of $L$ and is denoted by $\Ker(L)$.  It is always
contained in $\pr_1(L)$. When $\Ker(L)$ has constant fibre dimension,
it is integrable if and only if
\begin{equation}\label{eq:intcond}
\phi(X,Y,Z)=0 \;\;\; \mbox{ for all } X, Y \in \ker(\theta),\, Z \in \pr_1(L),
\end{equation}
at each point of $P$. In this case, the
leaves of the corresponding \textbf{characteristic foliation}
 are the null spaces of the
presymplectic forms along the leaves. On each leaf
$\iota:\mathcal{O}\hookrightarrow P$, the 2-form $\theta$
is basic with respect to the characteristic foliation if and only if
\begin{equation}\label{eq:basic}
\ker(\theta)\subseteq \ker(\iota^*\phi)
\end{equation}
at all points of $\mathcal{O}$. In this case,
forming the leaf space of this
foliation (locally, or globally when the foliation is simple) produces
a quotient manifold bearing a singular foliation by twisted symplectic leaves;
it is in fact a twisted Poisson manifold. In particular, when $\phi=0$,
conditions \eqref{eq:intcond} and \eqref{eq:basic} are satisfied,
and the quotient is an ordinary Poisson manifold.
Thus, Dirac manifolds can be
regarded as ``pre-Poisson'' manifolds, since, in nice situations, they
become Poisson manifolds after they are divided out by the
characteristic foliation.

Functions which are annihilated by all tangent vectors in the
characteristic distribution (equivalently, have differentials in the
projection of $L$ to $T^*P$)
are called \textbf{admissible} \cite{Cou90}.  For admissible
$f$ and $g$,
one can define
\begin{equation}\label{eq:admbracket}
\{f,g\}:=\theta(X_f,X_g),
\end{equation}
where $X_f$ is any vector field such that $(X_f,df) \in L$. (Note
that any two choices for $X_f$ differ by a characteristic vector,
so the bracket \eqref{eq:admbracket} is well defined.)
If \eqref{eq:basic} holds, then the algebra of
admissible functions is closed under this bracket, but it is not in
general a Poisson algebra, due to the presence of
$\phi$.
In particular, if the
characteristic foliation is regular and simple, the admissible
functions are just the functions on the (twisted) Poisson quotient.

\begin{example}(\textit{Nonintegrable characteristic distributions})

Consider the presymplectic structure $x_1 dx_1\wedge dx_2$ on
$\reals^2$.  Its characteristic distribution consists of the zero
subspace at points where $x_1 \neq 0$ and the entire tangent space at
each point of the $x_2$ axis.  Thus, the points off the axis are
integral manifolds, while there are no integral manifolds through
points on the axis.  The only admissible functions are constants.

On the other hand, if a 2-form is not closed, then its kernel
may have constant fibre dimension and still be nonintegrable.
For example, the characteristic distribution of the 2-form
$x_2dx_1\wedge dx_4 - dx_3\wedge dx_4$ on $\mathbb{R}^4$ is
spanned by ${\partial}/{\partial{x_1}}+x_2{\partial}/
{\partial{x_3}}$ and ${\partial}/{\partial{x_2}}$. A direct computation shows
that this 2-dimensional distribution does not satisfy the Frobenius
condition, so it is not integrable.
\end{example}

\begin{example}(\textit{A nonreducible 2-form})

The characteristic foliation of the 2-form $(x_3^2 + 1)dx_1\wedge
dx_2$ on $\mathbb{R}^3$ consists of lines parallel to the $x_3$-axis,
so it is simple.  However, the form is not basic with respect to
this foliation.
\end{example}

We will say more about presymplectic leaves and quotient Poisson structures in Section
\ref{subsec:diracmaps}.

\subsection{Poisson maps}\label{subsec:poissmap}

Although we shall see later that the following notion of morphism
between Poisson manifolds is not the only  one, it
is certainly the most obvious one.

Let $(P_1,\Pi_1)$ and $(P_2,\Pi_2)$ be Poisson manifolds. A smooth
map $\poismap:P_1 \to P_2$ is a \textbf{ Poisson map} if
$\poismap^*: C^\infty(P_2) \to C^\infty(P_1)$ is a homomorphism of
Poisson algebras, i.e.,
$$
\{f,g\}_2\circ \poismap = \{ f\circ \poismap, g \circ \poismap\}_1
$$
for $f,g \in C^\infty(P_2)$. One
can reformulate this condition in terms of Poisson bivectors or
hamiltonian vector fields as follows. A map $\poismap: P_1 \to P_2$ is a Poisson map if and only if
either of the following two equivalent conditions hold:
\begin{itemize}
\item[i)] $\poismap_* \Pi_1 = \Pi_2$, i.e., $\Pi_1$ and $\Pi_2$ are
$\poismap$-related.

\item[ii)] For all $f\in C^\infty(P_2)$, $X_f =
\poismap_*(X_{\poismap^*f})$.
\end{itemize}

It is clear by condition $ii)$ that  trajectories
of $X_{\poismap^*f}$ project to those of $X_f$ if $\poismap$ is a
Poisson map. This provides a way to ``lift'' some paths from $P_2$
to $P_1$. However, knowing that $X_f$ is complete does not
guarantee that $X_{\poismap^*f}$ is complete. In order to assure
that there are no ``missing points'' on the lifted trajectory on
$P_1$, we define a Poisson map $\poismap:P_1 \to P_2$ to be
\textbf{complete} if for any $f\in C^\infty(P_2)$ such that $X_f$
is complete, then $X_{\poismap^*f}$ is also complete.
Alternatively, one can replace the condition of $X_f$ being
complete by $X_f$ (or $f$ itself) having compact support.
Note that there is no notion of completeness (or
``missing point'') for a Poisson manifold by itself, only for a
Poisson manifold \textit{relative} to another.

\begin{remark} (\textit{Cotangent paths})

The path lifting alluded to above is best understood in terms of
so-called cotangent paths \cite{GG,We95}.  A \textbf{cotangent path} on a Poisson
manifold $P$ is a path $\alpha$ in $T^*P$ such that
$(\pi\circ\alpha)'= \widetilde{\Pi}\circ\alpha$, where $\pi$ is the
cotangent bundle projection.  If $\poismap:P_1 \to
P_2$ is a Poisson map, then a cotangent path $\alpha_1$ on $P_1$ is a
\textbf{horizontal lift} of the cotangent path $\alpha_2$ on $P_2$ if
$\alpha_1(t)=\poismap^*(\alpha_2(t))$ for all $t$.  It turns out that
a cotangent path on $P_2$ has at most one horizontal lift for each
initial value of $\pi\circ\alpha_1$.  Furthermore, the existence of
horizontal lifts
for all cotangent paths $\alpha_2$ and all initial data consistent
with the initial value of $\alpha_2$ is equivalent to completeness of
the map $\poismap$.

This path lifting property suggests that complete Poisson maps play
the role of ``coverings'' in Poisson geometry.  This idea is borne out
by some of the examples below.

\end{remark}

\begin{example} (\textit{Complete functions})

Let us regard $\mathbb{R}$ as a Poisson manifold, equipped with the
zero Poisson
bracket.
(This is the only possible Poisson structure on $\mathbb{R}$.)
Then any map $f:P \to \mathbb{R}$ is a Poisson map, which is
complete if and only if  $X_f$ is
a complete vector field.
\end{example}

Observe that the notion of completeness singles out the subset of $C^\infty(P)$ consisting of complete
functions, which is preserved under complete Poisson maps.

\begin{exercise}
For which Poisson manifolds is the set of complete functions closed
under addition?   (Hint: when are all functions complete?)
\end{exercise}

\begin{example}\label{ex:open}(\textit{Open subsets of symplectic manifolds})

Let $(P,\Pi)$ be a symplectic manifold, and let $U \subseteq P$ be an open subset.
Then the inclusion map $U \hookrightarrow P$ is complete if and only
if $U$ is closed (hence a union of connected components).  More
generally, the image of a complete Poisson map is a union of
symplectic leaves.
\end{example}

Example \ref{ex:open} suggests that (connected) symplectic manifolds are
``minimal objects'' among Poisson manifolds.

\begin{exercise}
The inclusion of every symplectic leaf in a Poisson manifold is a complete Poisson map.
\end{exercise}
\begin{exercise}
Let $P_1$ be a Poisson manifold, and let $P_2$ be symplectic. Then
any Poisson map $\poismap:P_1 \to P_2$ is a submersion.
Furthermore, if $P_2$ is connected and $\poismap$ is complete, then
$\poismap$ is surjective (assuming that $P_1$ is nonempty).
\end{exercise}

The previous exercise is the first step in establishing that
complete Poisson maps with symplectic target must be fibrations.
In fact, if $P_1$ is symplectic and $\mathrm{dim}(P_1) =
\mathrm{dim}(P_2)$, then  a complete Poisson map $\poismap:P_1 \to
P_2$ is a covering map. In general, a complete Poisson map
$\poismap:P_1 \to P_2$, where $P_2$ is \textit{symplectic}, is a
locally trivial symplectic fibration with a flat Ehresmann
connection: the horizontal lift in $T_x P_1$ of a vector $X$ in
$T_{\poismap(x)}P_2$ is defined as
$$
\widetilde{\Pi}_1((T_x\poismap)^*\widetilde{\Pi}_2^{-1}(X)).
$$
The horizontal subspaces define a foliation whose leaves
are coverings of $P_2$, and $P_1$ and $\poismap$ are completely determined, up
to isomorphism, by the holonomy
$$
\pi_1(P_2,x) \to \Aut(\poismap^{-1}(x)),
$$
see \cite[Sec.~7.6]{SilWein99} for details.


\subsection{Dirac maps}\label{subsec:diracmaps}
To see how to define Dirac maps, we first reformulate the
condition for a map $\poismap:(P_1,\Pi_1) \to (P_2,\Pi_2)$ to be
Poisson in terms of the bundles
$L_{\Pi_1}=\gra(\widetilde{\Pi}_1)$ and
$L_{\Pi_2}=\gra(\widetilde{\Pi}_2)$. First, note that
$\poismap$ is a Poisson map if and only if, at each $x\in P_1$,
\begin{equation}\label{eq:poissoncond}
(\widetilde{\Pi}_2)_{\poismap(x)}=T_x\poismap \circ
(\widetilde{\Pi}_1)_x \circ (T_x\poismap)^*.
\end{equation}
Now, using \eqref{eq:poissoncond}, it is not difficult to check
that $\poismap$ being a Poisson map is equivalent to
\begin{equation}\label{eq:poissbundlecond}
L_{\Pi_2}=\{(T\poismap(X),\beta) \, |\, X\in TP_1,\, \beta \in
T^*P_2, \, (X,(T\poismap)^*(\beta)) \in L_{\Pi_1} \}.
\end{equation}

Similarly, if $(P_1,\omega_1)$ and $(P_2,\omega_2)$ are
presymplectic manifolds, then a map $\poismap:P_1 \to P_2$
satisfies $\poismap^*\omega_2 = \omega_1$ if and only if
$L_{\omega_1}$ and $L_{\omega_2}$ are related by
\begin{equation}\label{eq:presymbundlecond}
L_{\omega_1}=\{ (X,(T\poismap)^*(\beta))\,|\, X \in TP_1,\,
\beta \in TP_2,\, (T\poismap(X),\beta) \in L_{\omega_2}\}.
\end{equation}

Since Dirac structures simultaneously generalize Poisson
structures and presymplectic forms, and  conditions
\eqref{eq:poissbundlecond} and \eqref{eq:presymbundlecond} both  make
sense for arbitrary Dirac subbundles, we have \textit{two} possible
definitions: If $(P_1,L_1)$ and $(P_2,L_2)$ are (possibly twisted)
Dirac manifolds, then a map $\poismap:P_1 \to P_2$ is a
\textbf{forward Dirac map} if
\begin{equation}\label{eq:forward}
L_2=\{(T\poismap(X),\beta) \, |\, X\in TP_1,\, \beta \in T^*P_2, \,
(X,(T\poismap)^*(\beta)) \in L_1 \},
\end{equation}
and a \textbf{backward Dirac map} if
\begin{equation}\label{eq:backward}
L_1=\{ (X,(T\poismap)^*(\beta))\,|\, X \in TP_1,\, \beta \in
TP_2,\, (T\poismap(X),\beta) \in L_2\}.
\end{equation}
Regarding vector Dirac structures as \textit{odd} (in the sense of
super geometry) analogues of lagrangian subspaces, one can
interpret formulas \eqref{eq:forward} and \eqref{eq:backward} via
composition of canonical relations \cite{We82}, see
\cite{BuRad02}.

The expression on the right-hand side of \eqref{eq:forward} defines
 at each point of $P_1$
a way to push a Dirac structure forward, whereas
\eqref{eq:backward} defines a pull-back operation. For this reason,
we often write
$$
L_2=\poismap_*L_1
$$
when \eqref{eq:forward} holds, following
the notation for $\poismap$-related
vector or bivector fields; similarly, we may write
$$
L_1=\poismap^*L_2
$$
instead of \eqref{eq:backward}.
This should explain the terminology ``forward'' and ``backward''.

\begin{remark}(\textit{Isotropic and coisotropic subspaces})

The notions of isotropic and coisotropic subspaces, as well as much of the
usual lagrangian\-/\- coisotropic calculus \cite{We82,We88} can be
naturally extended to Dirac vector spaces . This
yields an alternative characterization of forward (resp.
backward) Dirac maps in terms of their graphs being coisotropic
(resp. isotropic) subspaces of the suitable product Dirac space \cite{uchino}.
\end{remark}

Note that the pointwise pull back $\poismap^*L_2$ is always
a well-defined family of maximal isotropic subspaces in the fibres
of $TP_1\oplus T^*P_1$, though it may not be continuous,
whereas $\poismap_*L_1$ may not be well-defined at all.

\begin{exercise}
Consider a smooth map $f:P_1\to P_2$, and let
$L_2$ be a $\phi$-twisted Dirac structure on $P_2$.
Show that if $f^*L_2$ defines  a smooth vector bundle,
then its sections are automatically closed under the $f^*\phi$-twisted Courant bracket
on $P_1$ (so that $f^*L_2$ is a $f^*\phi$-twisted Dirac structure).
\end{exercise}

If $P_1$ and $P_2$ are symplectic manifolds, then a map $\poismap:
P_1 \to P_2$ is forward Dirac if and only if it is a Poisson map,
and backward Dirac if and only if it pulls back the symplectic
form on $P_2$ to the one on $P_1$, in which case we call it a
\textbf{symplectic map}.

The next example shows that
forward Dirac maps need not be backward Dirac, and vice versa.

\begin{example}\textit{(Forward vs. backward Dirac maps)}\label{ex:fversusb}

Consider $\mathbb{R}^2=\{(q,p)\}$, equipped with the symplectic form $dq\wedge dp$,
and $\mathbb{R}^4=\{(q_1,p_1,q_2,p_2)\}$, with symplectic form $dq_1\wedge dp_1 +
dq_2 \wedge dp_2$. Then a simple computation shows that
the inclusion
$$
\mathbb{R}^2 \hookrightarrow \mathbb{R}^4, \qquad (q,p)\mapsto
(q,p,0,0),
$$
is a symplectic (i.e. backward Dirac) map, but it does not preserve Poisson
brackets. On the other hand, the projection
$$
\mathbb{R}^4 \to \mathbb{R}^2, \;\;\; (q_1,p_1,q_2,p_2) \mapsto
(q_1,p_1),
$$
is a Poisson (i.e. forward Dirac) map, but it is not symplectic.
\end{example}

\begin{example}(\textit{Backward Dirac maps and
restrictions})\label{ex:backrest}

Let $(P,L)$ be a (possibly twisted) Dirac manifold, and let
$\iota:N \hookrightarrow P$ be a submanifold. Let $L_N \subset TN\oplus T^*N$
be the subbundle defined pointwise by the restriction of $L$ to $N$, see
\eqref{eq:canisom}, and suppose that $L_N$ is smooth, so that it defines
a Dirac structure on $N$. A direct computation shows that
$$
L_N = \iota^*L,
$$
hence the inclusion $\iota$ is a backward Dirac map.
\end{example}

The next exercise explains when the notions of forward and backward
Dirac maps coincide.

\begin{exercise}
Let $V_1$ and $V_2$ be vector spaces, and let $f:V_1 \to V_2$ be a linear map.
\begin{enumerate}
\item Let $L$ be a vector Dirac structure on $V_1$. Then $f^*f_*L
= L$ if and only if $\Ker(f)\subseteq \Ker(L)$, where
$\Ker(L)=V\cap L$.

\item Let $L$ be a vector Dirac structure on $V_2$. Then $f_*f^*L=L$
if and only if $f(V_1) \subseteq R$, where $R=\pr_1(L) \subseteq V_2$.
\end{enumerate}
It follows that $f^*f_*(L) = L$ for all $L$ if and only if $f$ is injective,
and $f_*f^*(L) = L$ for all $L$ if and only if $f$ is surjective.
\end{exercise}
In particular, the previous exercise shows that if $P_1$ and $P_2$
are symplectic manifolds, then a Poisson map $P_1\to P_2$
is symplectic if and only if it is an
immersion, and a symplectic map $P_1\to P_2$
is Poisson if and only if the map is a submersion (compare with
Example \ref{ex:fversusb}).  Thus, the only maps which are both
symplectic and Poisson are local diffeomorphisms.

Using the previous exercise, we find important examples of
maps which are both forward and backward Dirac.

\begin{example}(\textit{Inclusion of presymplectic leaves})\label{ex:inclus}

Let $(P,L)$ be a twisted Dirac manifold.  Let
$(\mathcal{O},\theta)$ be a presymplectic leaf, and let $\iota:
\mathcal{O} \hookrightarrow P$ be the inclusion. We regard
$\mathcal{O}$ as a Dirac manifold, with Dirac structure
$L_{\theta}=\gra(\widetilde{\theta})$. Then it follows from
the definition of $\theta$ that $\iota$ is a backward Dirac map.
On the other hand, since
$$
T\iota(T\mathcal{O}) = \pr_1(L)
$$
at each point, $\iota_*L_{\theta} = \iota_* \iota^*L=L$, so $\iota$ is also a forward Dirac map.

Note that $\theta$
 is completely determined by either of the
conditions that the inclusion be forward or backward Dirac.
\end{example}

\begin{example}(\textit{Quotient Poisson structures})\label{ex:quot}

Let $(P,L)$ be a Dirac manifold, and suppose that its characteristic
foliation is regular and simple.  According to  the discussion in
Section \ref{subsec:quotients}, the leaf space $P_{red}$ has an induced Poisson
structure $\Pi_{red}$. Using the definition of $\Pi_{red}$, one
can directly show that the natural projection
$$
\pr: P \longrightarrow P_{red}
$$
is a forward Dirac map, i.e., $\pr_*L =
\gra(\widetilde{\Pi}_{red})$. But since
$$
\Ker(T\pr) = \Ker(L),
$$
the previous exercise implies that $\pr^*\pr_*L=L$, so $\pr$ is a
backward Dirac map as well.

As in Example \ref{ex:inclus}, $\Pi_{red}$ is uniquely determined by
 either of the conditions that  $\pr$ be
backward or forward Dirac.
\end{example}

Example \ref{ex:quot} has an important particular case, which
illustrates the connection between Dirac geometry and
the theory of hamiltonian actions.

\begin{example}(\textit{Poisson reduction})\label{ex:reduction}

Suppose that $J:P \to \mathfrak{g}^*$ is the momentum map for a
hamiltonian action of a Lie group $G$ on a Poisson manifold $(P,\Pi)$.
Let $\mu \in \mathfrak{g}^*$ be a regular value for $J$, let
$Q=J^{-1}(\mu)$, and assume that the orbit space
$$
P_{red}=Q/G_{\mu}
$$
is a smooth manifold such that the projection $Q \to P_{red}$ is a
surjective submersion. Following Examples \ref{ex:momlevel} and
\ref{ex:backrest}, we know that $Q$ has an induced Dirac structure
$L_Q$ with respect to which the inclusion $ Q\hookrightarrow P$ is
a backward Dirac map.
\begin{exercise}
Show that the $G_{\mu}$-orbits on $Q$
coincide with the characteristic foliation of $L_Q$.
\end{exercise}

Thus, by Example
\ref{ex:quot}, $P_{red}$ inherits a Poisson structure $\Pi_{red}$
for which the projection $Q \to P_{red}$ is both backward and
forward Dirac (and either one of these conditions defines
$\Pi_{red}$ uniquely).
\end{example}

\section{Algebraic Morita equivalence}

There is another notion of morphism between Poisson manifolds which,
though it does {\em not} include all the Poisson maps, is
more closely adapted to the ``representation theory'' of Poisson
manifolds.   It is
based on an algebraic idea which we present first.  (The impatient
reader may skip to Section \ref{sec:geometric}.)

\subsection{Ring-theoretic Morita equivalence}\label{subsec:rings}

Let $\mathcal{A}$ and $\mathcal{B}$ be unital algebras over a
fixed ground ring $k$, and let $_{\st{\mathcal{A}}}\mathfrak{M}$
and $_{\st{\mathcal{B}}}\mathfrak{M}$ denote the categories of
left modules over $\mathcal{A}$ and $\mathcal{B}$, respectively.
We call $\mathcal{A}$ and $\mathcal{B}$ \textbf{Morita equivalent}
\cite{Morita58} if they have equivalent categories of left
modules, i.e., if there exist functors
\begin{equation}\label{eq:functors}
\mathcal{F}: {_{\st{\mathcal{B}}}\mathfrak{M}} \longrightarrow
{_{\st{\mathcal{A}}}\mathfrak{M}} \;\; \mbox{ and }\;\;
\widetilde{\mathcal{F}}:
{_{\st{\mathcal{A}}}\mathfrak{M}}\longrightarrow
{_{\st{\mathcal{B}}}\mathfrak{M}}
\end{equation}
whose compositions are naturally equivalent to the identity functors:
$$
{\mathcal{F}}\circ \widetilde{\mathcal{F}}\cong
\id_{_{\st{\mathcal{A}}}\mathfrak{M}},\;\;\mbox{ and }\;\;
\widetilde{\mathcal{F}} \circ {\mathcal{F}} \cong
\id_{_{\st{\mathcal{B}}}\mathfrak{M}}.
$$
One way to construct such functors between module categories
is via bimodules: if $\AXB$ is an
$(\mathcal{A},\mathcal{B})$-bimodule (i.e., $\mathrm{X}$ is a
$k$-module which is a left $\mathcal{A}$-module and a right
$\mathcal{B}$-module, and these actions commute), then we define an
associated functor
$\mathcal{F}_{\mathrm{X}}:{_{\st{\mathcal{B}}}\mathfrak{M}} \to
{_{\st{\mathcal{A}}}\mathfrak{M}}$ by setting, at the level of
objects,
\begin{equation}\label{eq:alginducton}
\mathcal{F}_{\mathrm{X}}({_{\st{\mathcal{B}}}\mathrm{M}}):= \AXB
\otimes_{\st{\mathcal{B}}}{_{\st{\mathcal{B}}}\mathrm{M}}.
\end{equation}
where the $\mathcal{A}$-module structure on
$\mathcal{F}_{\mathrm{X}}({_{\st{\mathcal{B}}}\mathrm{M}})$ is
 given by
$$
a\cdot (x \otimes_{\st{\mathcal{B}}} m)=
(ax)\otimes_{\st{\mathcal{B}}} m.
$$
For a morphism $T:{_{\st{\mathcal{B}}}\mathrm{M}} \longrightarrow
{_{\st{\mathcal{B}}}\mathrm{M}'}$, we define
\begin{equation}\label{eq:alginducton2}
\mathcal{F}_{\mathrm{X}}(T):\AXB
\otimes_{\st{\mathcal{B}}}{_{\st{\mathcal{B}}}\mathrm{M}}
\longrightarrow
\AXB\otimes_{\st{\mathcal{B}}}{_{\st{\mathcal{B}}}\mathrm{M}'},\;\;\;\;
\mathcal{F}_{\mathrm{X}}(T)(x\otimes_{\st{\mathcal{B}}}m)=
x\otimes_{\st{\mathcal{B}}}T(m).
\end{equation}

This way of producing functors turns out to be very general. In
fact, as we will see in Theorem  \ref{thm:morita}, any functor establishing an equivalence
between categories of modules is naturally equivalent to a functor
associated with a bimodule.

\begin{exercise}
Let $\mathrm{X}$ and $\mathrm{X}'$ be
$(\mathcal{A},\mathcal{B})$-bimodules. Show that the associated
functors $\mathcal{F}_{\mathrm{X}}$ and
$\mathcal{F}_{\mathrm{X}'}$ are   naturally equivalent
 if and only if the
bimodules $\mathrm{X}$ and $\mathrm{X}'$ are isomorphic.
\end{exercise}

It follows from the previous exercise that the functors $\mathcal{F}_{\X}:
{_{\st{\mathcal{B}}}\mathfrak{M}} \to
{_{\st{\mathcal{A}}}\mathfrak{M}}$, associated with an $(\A,\B)$-bimodule $\X$,
and $\mathcal{F}_{\Y}:
{_{\st{\mathcal{A}}}\mathfrak{M}} \to
{_{\st{\mathcal{B}}}\mathfrak{M}}$, associated with an $(\B,\A)$-bimodule $\Y$,
 are inverses of one another if and only if
\begin{equation}\label{eq:invbim}
\AXB \otimes_{\st{\mathcal{B}}} \BYA \cong \mathcal{A} \;\; \mbox{
and }\;\; \BYA \otimes_{\st{\mathcal{A}}}\AXB \cong \mathcal{B}.
\end{equation}
The isomorphisms in \eqref{eq:invbim} are \textit{bimodule}
isomorphisms, and $\mathcal{A}$ and $\mathcal{B}$ are regarded as
$(\mathcal{A},\mathcal{A})$- and
$(\mathcal{B},\mathcal{B})$-bimodules, respectively, in the
natural way (with respect to left and right multiplications).
So Morita equivalence is equivalent to the existence of bimodules
satisfying \eqref{eq:invbim}.

One can see Morita equivalence as the notion of isomorphism in an
appropriate category.  For that, we think of an \textit{arbitrary}
$(\mathcal{A},\mathcal{B})$-bimodule as a ``generalized morphism''
between $\mathcal{B}$ and $\mathcal{A}$. Note that, if
$\mathcal{A} \stackrel{q}{\leftarrow} \mathcal{B}$ is an ordinary
algebra homomorphism, then we can use it to make $\mathcal{A}$
into an $(\mathcal{A},\mathcal{B})$-bimodule by
\begin{equation}\label{eq:incl}
a \cdot x \cdot b := axq(b), \;\;\; a \in \mathcal{A}, \; x \in
\mathcal{A}, \; b \in \mathcal{B}.
\end{equation}
Since the tensor product
$$
\AXB \otimes_{\st{\mathcal{B}}} \BYC
$$
is an $(\A,\C)$-bimodule, we can see it as a  ``composition'' of
bimodules. As this composition is only associative up to
isomorphism, we consider the collection of \textit{isomorphism
classes} of $(\A,\B)$-bimodules, denoted by $\bim(\A,\B)$. Then
$\otimes_{\mathcal{B}}$ defines an associative composition
\begin{equation}\label{eq:composition}
\bim(\A,\B) \times \bim(\B,\C) \to \bim(\A,\C).
\end{equation}
We define the category $\catalg$ to be that
in which the objects are unital
$k$-algebras and the morphisms $\mathcal{A} \leftarrow
\mathcal{B}$ are the isomorphism classes of
$(\mathcal{A},\mathcal{B})$-bimodules, with composition given by
\eqref{eq:composition}; the identities are the algebras themselves
seen as  bimodules in the usual way. Note that a bimodule $\AXB$
is invertible in $\catalg$ if and only if it satisfies \eqref{eq:invbim}
for some bimodule $\BYA$, so the notion of
isomorphism in $\catalg$ coincides with Morita equivalence.

This is part of Morita's theorem \cite{Morita58},
see also \cite{Bass68}.
\begin{theorem}\label{thm:morita}
Let $\mathcal{A}$ and $\mathcal{B}$ be unital $k$-algebras.
\begin{enumerate}
\item A functor $\mathcal{F}: {_{\st{\mathcal{B}}}\mathfrak{M}}
\to {_{\st{\mathcal{A}}}\mathfrak{M}}$ is an equivalence of
categories if and only if there exists an invertible
$(\mathcal{A},\mathcal{B})$-bimodule $\mathrm{X}$ such that
$\mathcal{F}\cong \mathcal{F}_{\mathrm{X}}$.

\item A bimodule $\AXB$ is invertible if and only if it is
finitely generated and projective as a left $\mathcal{A}$-module
and as a right $\mathcal{B}$-module, and $\mathcal{A} \to
\End_{\mathcal{B}}(\mathrm{X})$ and $\mathcal{B} \to
\End_{\mathcal{A}}(\mathrm{X})$ are algebra isomorphisms.
\end{enumerate}
\end{theorem}

\begin{example}(\textit{Matrix algebras})\label{ex:matrix}

A unital algebra $\mathcal{A}$ is Morita equivalent to the matrix
algebra $M_n(\mathcal{A})$, for any $n\geq 1$, through the
$(M_n(\mathcal{A}),\mathcal{A})$-bimodule $\mathcal{A}^n$.
\end{example}

The following is a geometric example.

\begin{example}(\textit{Endomorphism bundles})\label{ex:endo}

Let $\A=C^\infty(M)$ be the algebra of complex-valued functions on
a manifold $M$. The Serre-Swan theorem asserts that any finitely
generated projective module over $C^\infty(M)$ can be identified with the
space of smooth sections $\Gamma(E)$ of a complex vector bundle $E \to M$.
In fact, $C^\infty(M)$ is Morita equivalent to $\Gamma(\mathrm{End}(E))$
via the $(\Gamma(\mathrm{End}(E)),C^\infty(M))$-bimodule $\Gamma(E)$.
When $E$ is the trivial bundle $\mathbb{C}^n\times M \to M$,
we recover the Morita equivalence of
$C^\infty(M)$ and $M_n(C^\infty(M))$ in Example \ref{ex:matrix}.
The same conclusion holds if $\mathcal{A}$ is the algebra of complex-valued
continuous functions on a compact Hausdorff space.
\end{example}

Morita equivalence preserves many algebraic properties
besides categories of representations, including ideal structures,
cohomology groups and deformation theories \cite{Bass68,GerSch}.
Another important Morita invariant is the center $\cent(\A)$ of a
unital
 algebra $\A$. If $\X$ is an invertible $(\A,\B)$-bimodule then,
for each $b\in \cent(\B)$, there is a unique $a=a(b) \in
\cent(\mathcal{A})$ determined by the condition $ax=xb$ for all $x
\in \X$. In this way, $\X$ defines an isomorphism
\begin{equation}\label{eq:projX}
\proj_{\X}: \cent(\A) \leftarrow \cent(\B),\;\;\; \proj_{\X}(b) =
a(b).
\end{equation}

The group of automorphisms of an object $\mathcal{A}$ in $\catalg$
is called its \textbf{Picard group}, denoted by
$\Pic(\mathcal{A})$. More generally, the invertible morphisms in
$\catalg$ form a ``large'' groupoid, called the \textbf{Picard
groupoid} \cite{Ben67}, denoted by $\Pic$. (Here, ``large'' refers to the fact
that the collection of objects in $\Pic$ is not a set,
though the collection of morphisms between any two of them is.)
The set of morphisms from $\mathcal{B}$ to $\mathcal{A}$ are the
Morita equivalences; we denote this set by
$\Pic(\mathcal{A},\mathcal{B})$. Of course
$\Pic(\mathcal{A},\mathcal{A})=\Pic(\mathcal{A})$. The orbit of an
object $\mathcal{A}$ in $\Pic$ is its Morita equivalence class,
while its isotropy $\Pic(\A)$ parametrizes the different ways
$\mathcal{A}$ can be Morita equivalent to any other object in its
orbit. It is clear from this picture that Picard groups of Morita
equivalent algebras are isomorphic.

Let us investigate the difference between $\Aut(\mathcal{A})$, the
group of ordinary algebra automorphisms of $\A$, and
$\Pic(\mathcal{A})$. Since ordinary automorphisms of $\A$ can be
seen as generalized ones, see \eqref{eq:incl}, we obtain a group
homomorphism
\begin{equation}\label{eq:jmap}
\inc: \Aut(\A) \to \Pic(\A).
\end{equation}
A simple computation shows that $\ker(\inc)=\Inaut(\A)$, the group
of inner automorphisms of $\A$. So the outer automorphisms
$\Outaut(\A):=\Aut(\A)/\Inaut(\A)$ sit inside $\Pic(\A)$.

\begin{exercise}
Morita equivalent algebras have isomorphic Picard groups. Do they
always have isomorphic groups of outer automorphisms?  (Hint: consider
the direct sum of two matrix algebras of the same or different sizes.)
\end{exercise}

On the other hand, \eqref{eq:projX} induces a group homomorphism
\begin{equation}\label{eq:proj}
\proj: \Pic(\A) \to \Aut(\cent(\A)),
\end{equation}
whose kernel is denoted by $\SPic(\A)$, the \textbf{static Picard
group} of $\A$.

\begin{remark}
If $\A$ is commutative, then each invertible bimodule induces an
automorphism of $\A$ by \eqref{eq:projX}, and $\SPic(\A)$ consists
of those bimodules ``fixing'' $\A$, which motivates our
terminology. Bimodules in $\SPic(\A)$ can also be characterized by
having equal left and
right module structures, and $\SPic(\A)$ is often
referred to in the literature as the ``commutative'' Picard group
of $\A$.
\end{remark}
If $\A$ is commutative, then the composition
$$
\Aut(\A) \stackrel{\inc}{\to} \Pic(\A) \stackrel{\proj}{\to}
\Aut(\A)
$$
is the identity. As a result, we can write $\Pic(\A)$ as a semi-direct product,
\begin{equation}\label{eq:piccomm}
\Pic(\A) = \Aut(\A) \ltimes \SPic(\A).
\end{equation}
The action of $\Aut(\A)$ on $\SPic(\A)$ is given by $\X
\stackrel{q}{\mapsto} {{_{\st{q}}}\X_{\st{q}}}$, where the left
and right $\A$-module structures on ${{_{\st{q}}}\X_{\st{q}}}$ are
$a\cdot x := q(a)x$ and $x \cdot b := x q(b)$. Although the orbits
of \textit{commutative} algebras in $\Pic$ are just their
isomorphism classes in the ordinary sense, \eqref{eq:piccomm}
illustrates that their isotropy groups in $\Pic$ may be
bigger than their ordinary automorphism groups. The following is a
geometric example.

\begin{example}(\textit{Picard groups of algebras of functions})\label{ex:picfn}

Let $\A = C^\infty(M)$ be the algebra of smooth
complex-valued functions on a manifold $M$. Using the Serre-Swan
identification of smooth complex
vector bundles over $M$ with projective modules
over $\A$, one can check that $\SPic(\A)$ coincides with
$\Pic(M)$, the group of isomorphism classes of complex line
bundles on $M$, which is isomorphic to $H^2(M,\mathbb{Z})$ via
the Chern class map. We then have a purely geometric description
of $\Pic(\A)$ as
\begin{equation}\label{eq:picfunc}
\Pic(C^\infty(M)) = \Diff(M) \ltimes
H^2(M,\mathbb{Z}),
\end{equation}
where the action of $\Diff(M)$ on $H^2(M,\mathbb{Z})$ is given by
pull back. In \eqref{eq:picfunc}, we use the identification of
algebra automorphisms of $\mathcal{A}$ with diffeomorphisms of
$M$, see e.g. \cite{Mrcun}.
\end{example}

\subsection{Strong Morita equivalence of $C^*$-algebras}

The notion of Morita equivalence of unital algebras has been
adapted to several other classes of algebras. An example is the notion of
\textit{strong Morita equivalence} of $C^*$-algebras, introduced by
Rieffel in \cite{Ri74a,Ri74}.

A \textbf{$C^*$-algebra} $\A$ is a complex
Banach algebra with an involution $^*$ such that
$$
\| a a^*\| = \|a\|^2, \;\; a \in \A.
$$
Important examples are the algebra of complex-valued continuous
functions on a locally compact Hausdorff space  and $\mathcal{B}(\mathcal{H})$,
the algebra of bounded operators on a Hilbert space $\mathcal{H}$.

The relevant category of modules over a $C^*$-algebra, to be
preserved under strong Morita equivalence, is that  of Hilbert
spaces on which the $C^*$-algebra acts through bounded operators.
More precisely, for a given $C^*$-algebra $\mathcal{A}$, we
consider the category $\repC(\mathcal{\A})$ whose objects are
pairs $(\mathcal{H},\rho)$, where $\mathcal{H}$ is a Hilbert space
and $\rho:\mathcal{A} \to \mathfrak{B}(\mathcal{H})$ is a
\textit{nondegenerate} $^*$-homomorphism of algebras, and morphisms are
bounded linear intertwiners. (Here ``nondegenerate'' means that
$\rho(\A)h=0$ implies that $h=0$, which is always satisfied if $\A$ is unital
and $\rho$ preserves the unit.)

Since we are now dealing with more elaborate modules,
it is natural that a bimodule giving rise to a functor
$\repC(\mathcal{B}) \to \repC(\A)$ analogous to
\eqref{eq:alginducton} should be equipped with extra
structure. If $(\mathcal{H},\rho) \in \repC(\B)$ and $\AXB$ is an
$(\A,\B)$-bimodule, the key observation is that if $\X$ is itself
equipped with an inner product $\SP{\cdot,\cdot}_{\st{\B}}$
\textit{with values in} $\B$, then the map $\AXB \otimes_{\st{\B}}
\mathcal{H} \times \AXB \otimes_{\st{\B}}\mathcal{H} \to \mathbb{C}$
uniquely defined by
\begin{equation}
(x_1\otimes h_1,x_2\otimes h_2) \mapsto
\SP{h_1,\rho(\SP{x_1,x_2}_{\st{\B}})h_2}
\end{equation}
is an inner product on $\AXB \otimes_{\st{\B}}
\mathcal{H}$, which we can complete to obtain  a Hilbert
space $\mathcal{H}'$. Moreover, the natural $\A$-action on $\AXB
\otimes_{\st{\B}} \mathcal{H}$ gives rise to a $^*$-representation
$\rho': \A \to \mathfrak{B}(\mathcal{H}')$. These are the main ingredients of
Rieffel's induction of representations \cite{Ri74a}.

More precisely, let $\X$ be a right $\B$-module. Then a $\B$-valued
inner product $\SP{\cdot,\cdot}_{\st{\B}}$ on $\X$ is a
$\mathbb{C}$-sesquilinear pairing $\X \times \X \to \B$ (linear in
the second argument) such that, for all $x_1,x_2 \in \X$ and $b
\in \B$, we have
$$
\SP{x_1,x_2}_{\st{\B}} = \SP{x_2,x_1}_{\st{\B}}^*, \;\;
\SP{x_1,x_2b}_{\st{\B}}=\SP{x_1,x_2}_{\st{\B}}b,\; \mbox{ and }
\SP{x_1,x_1}_{\st{\B}} > 0\; \mbox{ if } x_1 \neq 0.
$$
(Inner products on left modules are defined analogously, but
linearity is required in the first argument). One can show that
$\|x\|_{\st{\B}}:=\|\SP{x,x}_{\st{\B}}\|^{1/2}$ is a norm in $\X$.
A (right) \textbf{Hilbert $\B$-module} is a (right) $\B$-module
$\X$ together with a $\B$-valued inner product
$\SP{\cdot,\cdot}_{\st{\B}}$ so that $\X$ is complete with respect
to $\|\cdot\|_{\st{\B}}$. Just as for Hilbert spaces, we denote by
$\mathfrak{B}_{\st{\B}}(\X)$ the algebra of endomorphisms of $\X$
possessing an adjoint with respect to
$\SP{\cdot,\cdot}_{\st{\B}}$.

\begin{example}(\textit{Hilbert spaces})\label{ex:hilbert}

If $\B=\mathbb{C}$, then Hilbert $\B$-modules are just ordinary
Hilbert spaces. In this case, $\mathfrak{B}_{\mathbb{C}}(\X)$
coincides with the algebra of bounded linear operators on $\X$, see e.g. \cite{RaWi}.
\end{example}

\begin{example}(\textit{Hermitian vector bundles})

Suppose $\B=C(X)$, the algebra of complex-valued continuous
functions on a compact Hausdorff space $X$. If
$E\to X$ is a complex vector bundle equipped with a hermitian metric $h$,
then $\Gamma(E)$ is a
Hilbert $\B$-module with respect to the $C(X)$-valued inner product
$$
\SP{e,f}_{\st{\B}}(x):=h_x(e(x),f(x)).
$$

To describe the most general Hilbert modules over $C(X)$, one needs
Hilbert bundles, which recover Example \ref{ex:hilbert} when $X$ is a point.
\end{example}

\begin{example}\label{ex:cstar}(\textit{$C^*$-algebras})

Any $C^*$-algebra $\B$ is a Hilbert $\B$-module with respect to the inner
product $\SP{b_1,b_2}_{\st{\B}}=b_1^*b_2$.
\end{example}

As in  the case of unital algebras, one can define,
for $C^*$-algebras $\A$ and $\B$, a ``generalized morphism''
$\A \leftarrow \B$ as a right Hilbert $\B$-module
$\X$, with inner product $\SP{\cdot,\cdot}_{\st{\B}}$, together
with a \textit{nondegenerate} $^*$-homomorphism $\A \to
\mathfrak{B}_{\st{\B}}(\X)$.
We ``compose'' $\AXB$ and $\BYC$ through a more elaborate tensor product:
we consider the algebraic
tensor product $\AXB \otimes_{\mathbb{C}} \BYC$, equipped with the
semi-positive $\mathcal{C}$-valued inner product uniquely defined by
\begin{equation}\label{eq:innerprod}
(x_1\otimes y_1,x_2\otimes y_2) \mapsto \SP{y_1,
\SP{x_1,x_2}_{\st{\mathcal{B}}}y_2}_{\st{\mathcal{C}}}.
\end{equation}
The null space of this inner product coincides with the span of
elements of the form $xb\otimes y- x\otimes by$ \cite{Lance}, so
\eqref{eq:innerprod} induces a positive-definite
$\mathcal{C}$-valued inner product on $\AXB
\otimes_{\st{\mathcal{B}}}\BYC$. The completion
of this space with respect to $\|\cdot\|_{\st{\mathcal{C}}}$ yields
a ``generalized morphism'' from $\mathcal{C}$ to $\A$
denoted by $\AXB\widehat{\otimes}_{\st{\B}} \BYC$,
called the \textbf{Rieffel tensor product} of $\AXB$ and $\BYC$.

An isomorphism between ``generalized morphisms'' is a bimodule isomorphism
preserving inner products. Just as ordinary tensor
products, Rieffel tensor products are
associative up to natural isomorphisms. So one can define a
category $\catCstar$ whose objects are
$C^*$-algebras and whose morphisms are isomorphism classes of
``generalized morphisms'', with composition given by Rieffel tensor
product; the identities are the algebras themselves, regarded as
bimodules in the usual way, and with the inner product of Example
\ref{ex:cstar}.

 Two $C^*$-algebras are \textbf{strongly Morita
equivalent} if they are isomorphic in $\catCstar$. As in the case
of unital algebras, isomorphic $C^*$-algebras are necessarily strongly
Morita equivalent.

\begin{remark}(\textit{Equivalence bimodules})

The definition of strong Morita equivalence
as isomorphism in $\catCstar$ coincides with Rieffel's original
definition in terms of equivalence bimodules (also called imprimitivity bimodules)
\cite{Ri74a,Ri74}. In fact, any ``generalized morphism'' $\AXB$ which is invertible
in $\catCstar$ can be endowed with an $\A$-valued inner
product, compatible with its $\B$-valued inner product in the appropriate way,
 making it into an equivalence bimodule, see \cite{Land01} and references therein.
Conversely, any equivalence bimodule
is automatically invertible in $\catCstar$.
\end{remark}

\begin{example}(\textit{Compact operators})

A Hilbert space $\mathcal{H}$, seen as a bimodule for $\mathbb{C}$ and
 the $C^*$-algebra $\mathcal{K}(\mathcal{H})$ of compact
operators on $\mathcal{H}$, defines a strong Morita equivalence.
\end{example}

\begin{example}(\textit{Endomorphism bundles})

Analogously to Example \ref{ex:endo}, a hermitian vector bundle $E\to X$,
where $X$ is a compact Hausdorff space, defines a strong Morita equivalence
between $\Gamma(\mathrm{End}(E))$ and $C(X)$.
\end{example}

Any ``generalized morphism''  $\AXB$ in $\catCstar$ defines a functor
$$
\mathcal{F}_{\X}:\repC(\B) \to \repC(\A),
$$
similar to \eqref{eq:alginducton}, but with Rieffel's tensor
product replacing the ordinary one, i.e., on objects,
\begin{equation}\label{eq:rieffelinduction}
\mathcal{F}_{\X}(\mathcal{H}):= \AXB \widehat{\otimes}_{\st{\B}}
\mathcal{H}.
\end{equation}
Such a functor is called \textbf{Rieffel induction of
representations} \cite{Ri74a}. It follows that strongly Morita
equivalent $C^*$-algebras have equivalent categories of
representations, although, in this setting, the converse is not
true \cite{Ri74} (see \cite{Blecher} for a different approach
where a converse does hold).

\begin{remark}(\textit{Strong vs. ring-theoretic Morita equivalence})

By regarding unital $C^*$-algebras simply as unital algebras over
$\mathbb{C}$, one can compare strong and ring-theoretic Morita
equivalences. It turns out that two unital $C^*$-algebras are
strongly Morita equivalent if and only if they are Morita
equivalent as unital $\mathbb{C}$-algebras \cite{Beer}. However, the Picard
groups associated to each notion are different in general, see
\cite{BuWa03}. In terms of Picard groupoids, this means that, over
unital $C^*$-algebras, the Picard groupoids associated with
ring-theoretic and strong Morita equivalences have the same
orbits, but generally different isotropy groups.
\end{remark}

A study of Picard groups associated with strong
Morita equivalence, analogous to
the discussion in Section \ref{subsec:rings}, can be found in \cite{BGR}.

\subsection{Morita equivalence of deformed algebras}\label{subsec:deformed}

Let $(P,\Pi)$ be a Poisson manifold and $C^\infty(P)$
be its algebra of smooth complex-valued functions.
The general idea of a \textbf{deformation quantization} of $P$ ``in
the direction'' of $\Pi$ is that of a family
$\star_{\hbar}$ of associative algebra structures on
$C^\infty(P)$ satisfying the following two conditions:
\begin{itemize}
\item[$i.)$] $f \star_{\hbar} g = f\cdot g + O(\hbar)$;

\item[$ii.)$]$\frac{1}{i\hbar}(f \star_{\hbar} g - g \star_{\hbar} f)
\longrightarrow
 \{f,g\},$ when $\hbar \to 0$.
\end{itemize}
There are several versions of deformation quantization. We will consider

\begin{enumerate}
\item {\bf Formal deformation quantization} \cite{BFFLS}: In this case,
$\star_{\hbar}$ is an associative product on
$C^\infty(P)[[\hbar]]$, the space of formal power series with coefficients
in $C^\infty(P)$.  Here $\hbar$ is a formal parameter, and
the ``limit'' in $ii.)$ above is defined simply by setting $\hbar$ to $0$.
A formal deformation quantization is also called a \textbf{star product}.
The contribution by Cattaneo and Indelicato \cite{CaIn} to
this volume contains a thorough exposition of the theory of star
products and its history.

\item {\bf Rieffel's strict deformation quantization} \cite{Ri89}:
In this setting, one starts with a dense Poisson subalgebra of
$C_\infty(P)$, the $C^*$-algebra of continuous functions on $P$
vanishing at infinity, and considers families of associative
products $\star_{\hbar}$ on it, defined along with norms and
involutions such that the completions form a continuous field of
$C^*$-algebras. The parameter $\hbar$ belongs to a closed subset
of $\mathbb{R}$ having $0$ as a non-isolated point, and  one can
make analytical sense of the limit in $ii.)$ above. Variations of
Rieffel's notion of deformation quantization are discussed in
\cite{Land98}.
\end{enumerate}

Intuitively, one should regard a deformation quantization $\star_{\hbar}$ as a
path in the ``space of associative algebra structures'' on $C^\infty(P)$
for which the Poisson structure $\Pi$ is the ``tangent vector'' at $\hbar = 0$.
 From this perspective, a direct relationship
between deformation quantization and Poisson geometry is more likely
in the formal case.

A natural question is when two algebras obtained by deformation quantization  are
Morita equivalent.
In the framework of formal deformation quantization, the first
observation is that if two deformation quantizations
$(C^\infty(P_1,\Pi_1)[[\hbar]],\star_{\hbar}^1)$ and
$(C^\infty(P_2,\Pi_2)[[\hbar]],\star_{\hbar}^2)$ are Morita
equivalent (as unital algebras over  $\mathbb{C}[[\hbar]]$), then the
underlying Poisson manifolds are isomorphic. So we can restrict
ourselves to a fixed Poisson manifold.
The following result is proven in
 \cite{Bu,BuWa2002}:

\begin{theorem}\label{thm:BuWa}
Let $P$ be symplectic. If $\Pic(P) \cong H^2(P,\mathbb{Z})$ has no torsion, then
it acts freely on the set of equivalence classes of star
products on $P$, and two star products are Morita equivalent if and only
if their classes lie in the same $H^2(P,\mathbb{Z})$-orbit,
up to symplectomorphism.
\end{theorem}

Recall that two star products  $\star_{\hbar}^1$ and $\star_{\hbar}^2$
are {\bf equivalent} if there exists a family of differential operators
$T_r:C^\infty(P)\to C^\infty(P)$, $r=1,2 \ldots$, so that
$T=\id + \sum_{r=1}^\infty T_r \hbar^r$ is an algebra isomorphism
$$
(C^\infty(P)[[\hbar]],\star_{\hbar}^1) \stackrel{\sim}{\longrightarrow}
(C^\infty(P)[[\hbar]],\star_{\hbar}^2).
$$
Equivalence classes of star products on a symplectic manifold
are parametrized by elements in
\begin{equation}\label{eq:chclass}
\frac{1}{i\hbar}[\omega]+ H^2_{\st{dR}}(P)[[\hbar]],
\end{equation}
where $\omega$ is the symplectic form on $P$ and
$H^2_{\st{dR}}(P)$ is the second de Rham cohomology group of $P$
with complex coefficients \cite{BFFLS,fe:simple,ne-ts:algebraic},
called {\bf characteristic classes}. As shown in \cite{BuWa2002},
the $\Pic(P)$-action of Theorem \ref{thm:BuWa} is explicitly given
in terms of these  classes by
\begin{equation}
[\omega_{\hbar}] \mapsto [\omega_{\hbar}] + 2\pi i c_1(L),
\end{equation}
where $[\omega_{\hbar}]$ is an element in \eqref{eq:chclass} and $c_1(L)$
is the image of the Chern class of the line bundle $L$ in $H^2_{\st{dR}}(P)$.

\begin{remark}\label{rem:formality}
A version of Theorem \ref{thm:BuWa} holds for arbitrary
Poisson manifolds $(P,\Pi)$, see \cite{Bu,JSW2001}. In this general
setting, equivalence classes of star products are parametrized by
classes of formal Poisson bivectors $\Pi_{\hbar}= \Pi + \hbar\Pi_1 +
\cdots$ (see \cite{Kont} or the exposition in \cite{CaIn}),
and the $\Pic(P)$-action on them, classifying
Morita equivalent deformation quantizations of $P$, is via gauge
transformations (see Section \ref{subsec:gauge}).
\end{remark}


In the framework of strict deformation quantization and the special
case of tori, a classification result for Morita equivalence was
obtained by Rieffel and Schwarz in \cite{RiefSch} (see also
\cite{Li,WeTan}).  Let us consider
$\mathbb{T}^n=\mathbb{R}^n/\mathbb{Z}^n$ equipped with a constant
Poisson structure, represented by a skew-symmetric real matrix $\Pi$:
if $(\theta_1,\ldots,\theta_n)$ are coordinates on ${T}^n$,
then
$$
\Pi_{ij}=\{\theta_i,\theta_j\}.
$$

Via the
Fourier transform, one can identify the algebra $C^\infty(\mathbb{T}^n)$
with  the space $\mathcal{S}(\mathbb{Z}^n)$  of complex-valued
functions on $\mathbb{Z}^n$ with rapid decay at infinity.  Under
 this identification, the pointwise
product of functions becomes the convolution on
$\mathcal{S}(\mathbb{Z}^n)$,
$$
\widehat{f}*\widehat{g}(n)= \sum_{k\in
\mathbb{Z}^n}\widehat{f}(n)\widehat{g}(n-k),
$$
$\widehat{f},\widehat{g} \in \mathcal{S}(\mathbb{Z}^n)$. One can
now use the matrix $\Pi$ to ``twist'' the convolution and define a
new product
\begin{equation}\label{eq:twistedconv}
\widehat{f}*_{\hbar}\widehat{g}(n)= \sum_{k\in
\mathbb{Z}^n}\widehat{f}(n)\widehat{g}(n-k)e^{-\pi i\hbar
\Pi(k,n-k)}
\end{equation}
on $\mathcal{S}(\mathbb{Z}^n)$, which can be pulled back to a new product in
$C^\infty(\mathbb{T}^n)$.  Here $\hbar$ is a real parameter. If we
set $\hbar =1$, this defines the algebra $\A_{\Pi}^\infty$, which can be thought
of as the ``algebra of smooth functions on the quantum torus
$\mathbb{T}^n_{\Pi}$''.
A suitable completion of $\A_{\Pi}^\infty$ defines a
$C^*$-algebra $\mathcal{A}_{\Pi}$, which is then thought of as the
``algebra of continuous functions on
$\mathbb{T}^n_{\Pi}$''.  (Note that, with $\hbar=1$, we are no longer
really considering a deformation.)

\begin{exercise}
Show that $1$ is a unit for $\A_{\Pi}$. Let $u_j = e^{2\pi i
\theta_j}$. Show that $u_j *_1 \bar{u}_j=\bar{u}_j *_1 u_j = 1$
and
\begin{equation}\label{eq:relation}
u_j *_1 u_k = e^{2\pi i \Pi_{jk}} u_k *_1 u_j.
\end{equation}
\end{exercise}
The algebra $\A_{\Pi}$ can be alternatively described as the universal
$C^*$-algebra generated by $n$ unitary elements $u_1,\ldots,u_n$
subject to the commutation relations \eqref{eq:relation}.

In this context, the question to be addressed is when
skew-symmetric matrices $\Pi$ and $\Pi'$ correspond to Morita
equivalent $C^*$-algebras $\A_{\Pi}$ and $\A_{\Pi'}$. Let
$O(n,n|\mathbb{R})$ be the group of linear automorphisms of
$\mathbb{R}^n\oplus {\mathbb{R}^n}^*$ preserving the inner product
\eqref{eq:symbr}. One can identify elements of $O(n,n|\mathbb{R})$
with matrices
$$
g=\left(%
\begin{array}{cc}
  A & B \\
  C & D \\
\end{array}%
\right),
$$
where $A,B,C$ and $D$ are $n\times n$ matrices satisfying
$$
A^tC + C^tA = 0 = B^tD + D^tB,\; \mbox{ and }\;  A^tD + C^tB =1.
$$
The group $O(n,n|\mathbb{R})$ ``acts'' on the space of
all $n\times n$ skew-symmetric matrices by
\begin{equation}\label{eq:fractional}
\Pi \;  {\mapsto}\;  g\cdot \Pi := (A\Pi + B)(C\Pi + D)^{-1}.
\end{equation}
Note that this is not an honest action, since the formula above
only makes sense when $(C\Pi + D)$ is invertible.

Let $SO(n,n|\mathbb{Z})$ be the subgroup of $O(n,n|\mathbb{R})$
consisting of matrices with integer coefficients and determinant
1. The main result of \cite{RiefSch}, as improved in
 \cite{Li,WeTan}, is

\begin{theorem}\label{thm:MEqtori}
If $\Pi$ is a skew-symmetric matrix, $g \in SO(n,n|\mathbb{Z})$ and $g\cdot \Pi$ is defined,
then $\A_{\Pi}$ and $\A_{g\cdot\Pi}$ are strongly Morita equivalent.
\end{theorem}

\begin{remark}(\textit{Converse results})

The converse of Theorem \ref{thm:MEqtori} holds for $n=2$ \cite{Ri81},
but not in general. In fact, for $n=3$, one can find $\Pi$ and $\Pi'$,
\textit{not} in the same $SO(n,n|\mathbb{Z})$-orbit,
for which $\A_{\Pi}$ and $\A_{\Pi'}$ are isomorphic (hence Morita equivalent)
\cite{RiefSch}.

On the other hand, for \textit{smooth} quantum tori, Theorem
\ref{thm:MEqtori} and its converse
hold with respect to a refined
notion of Morita equivalence, called
``complete Morita equivalence'' \cite{Sch}, in which bimodules
carry connections of constant curvature.

For the algebraic Morita
equivalence of smooth quantum tori, see \cite{ElLi}.
\end{remark}

\begin{remark}(\textit{Dirac structures and quantum tori})

In \cite{RiefSch}, the original version of Theorem \ref{thm:MEqtori}
was proven under an additional hypothesis. Rieffel and Schwarz consider
three types of generators of $SO(n,n|\mathbb{Z})$, and prove that
their action preserves Morita equivalence. In order to show that
$\A_{\Pi}$ and $\A_{g\Pi}$ are Morita equivalent for an arbitrary
$g \in SO(n,n|\mathbb{Z})$ (for which $g\Pi$ is defined), they need
to assume that $g$ can be written as a product of generators $g_r\cdots g_1$
in such a way that each of the products $g_k\cdots g_1\Pi$ is defined.
The result in Theorem \ref{thm:MEqtori}, without this assumption, is conjectured
in \cite{RiefSch}, and it was proven by Li in \cite{Li}.

A geometric way to circumvent the difficulties in the Rieffel-Schwarz proof,
in which Dirac structures play a central role, appears in \cite{WeTan}.
The key point is the observation that, even if $g\cdot\Pi$ is not
defined as a skew-symmetric matrix,
it is still a Dirac structure on $\mathbb{T}^n$.
The authors develop a way to quantize constant Dirac structures on $\mathbb{T}^n$
by
attaching to each one of them a
Morita equivalence class of quantum tori.
They extend the $SO(n,n|\mathbb{Z})$ action
to Dirac structures and prove that the
Morita equivalence classes of the corresponding quantum tori is unchanged under
the action.
\end{remark}

\section{Geometric Morita equivalence}\label{sec:geometric}

In this section, we introduce a purely geometric notion of Morita
equivalence of Poisson manifolds.  This notion leads inevitably to
the consideration of Morita equivalence of symplectic groupoids, so we
will make a digression into the Morita theory of general Lie groups
and groupoids.  We end the section with a discussion of gauge
equivalence, a geometric equivalence which is close to Morita
equivalence, but is also related to the algebraic Morita equivalence
of star products, as discussed in Section \ref{subsec:deformed}.

\subsection{Representations and tensor product}\label{subsec:geomrep}
In order to define Morita equivalence in Poisson
geometry,  we need notions of ``representations'' of (or
``modules'' over) Poisson manifolds as well as their tensor products.

As we saw in Example \ref{ex:open}, symplectic manifolds are in some sense
``irreducible'' among Poisson manifolds. If one thinks of
Poisson manifolds as algebras, then symplectic manifolds could
be thought of as ``matrix algebras''.
Following this analogy, a representation of a Poisson manifold $P$ should be  a
symplectic manifold $S$ together with a Poisson map $J:S \to P$
which is complete. At the level of functions, we have a
``representation'' of $C^\infty(P)$ by $J^*:C^\infty(P) \to
C^\infty(S)$. This notion of representation is also suggested
by the theory of geometric quantization, in which symplectic manifolds
become ``vector spaces'' on which their Poisson algebras ``act
asymptotically''.

More precisely, we define a left [right] \textbf{$P$-module} to be
 a complete [anti-]
symplectic realization $J:S \to P$. Our first example illustrates how
modules over Lie-Poisson manifolds are related to hamiltonian actions.

\begin{example}(\textit{Modules over $\mathfrak{g}^*$ and hamiltonian
actions})\label{ex:gstarmod}

Let $(S,\Pi_S)$ be a symplectic Poisson manifold, $\mathfrak{g}$ be a
Lie algebra, and suppose that
$J:S \to \mathfrak{g}^*$ is a symplectic realization of $\mathfrak{g^*}$.
The map
\begin{equation}
\mathfrak{g} \to \mathcal{X}(S),\;\;\; v \mapsto
\widetilde{\Pi}_S(dJ_v),
\end{equation}
where $J_v(x):=\SP{J(x),v}$, defines a $\mathfrak{g}$-action on
$S$ by hamiltonian vector fields for which $J$ is the momentum
map. On the other hand, the momentum map $J:S\to \mathfrak{g}^*$
for a hamiltonian $\mathfrak{g}$-action on $S$ is a Poisson map,
so we have a one-to-one correspondence between symplectic
realizations of $\mathfrak{g}^*$ and hamiltonian
$\mathfrak{g}$-manifolds.

A symplectic realization $J:S\to \mathfrak{g}^*$ is complete if and only
if the associated infinitesimal hamiltonian action is by complete vector
fields, in which case it can be integrated to a hamiltonian $G$-action,
where $G$ is the connected and simply-connected Lie group
 having $\mathfrak{g}$
as its Lie algebra.
So $\gstar$-modules are just the same thing as hamiltonian $G$-manifolds.
\end{example}

\begin{remark}(\textit{More general modules over $\mathfrak{g}^*$})\label{rem:genmod}

The one-to-one correspondence in Example \ref{ex:gstarmod} extends
to one
between Poisson maps into $\mathfrak{g}^*$ (from any Poisson manifold, not necessarily
symplectic) and hamiltonian $\mathfrak{g}$-actions on Poisson manifolds, or, similarly,
between complete Poisson maps into $\mathfrak{g}^*$ and Poisson manifolds carrying
hamiltonian $G$-actions.
This indicates that it may be useful to regard \textit{arbitrary} (complete) Poisson maps as
modules over Poisson manifolds; we will say more about this in Remarks
\ref{rem:hamiltpoiss} and \ref{rem:general}.
\end{remark}

We now define a tensor product operation on modules over a Poisson manifold.
Let $J:S \to P$ be a right $P$-module, and let $J':S'\to P$
be a left $P$-module.  Just as, in algebra, we can think of the tensor
product over $\A$ of a left module $\X$ and a right module $\Y$ as a
quotient of their tensor product over the ground ring $k$, so in
Poisson geometry we can define the tensor product of $S$ and $S'$ to
be a ``symplectic quotient'' of $S \times S'$.  Namely, the fibre
product
\begin{equation}\label{eq:fibreprod}
S \times_{(J,J')} S' = \{(x,y) \in S \times S'
\;|\; J(x)=J'(y)\}
\end{equation}
is the inverse image of the diagonal under the Poisson map
$(J,J'):S\times S'\to P\times \overline{P}$, hence, whenever it is
smooth, it is a coisotropic submanifold. (Here $\overline{P}$
denotes $P$ equipped with its Poisson structure multiplied by -1.)
Let us assume then, that the fibre product is smooth; this is the
case, for example, if either $J$ or $J'$ is a surjective
submersion.  Then we may define the \textbf{tensor product} $S*S'$
over $P$ to be the quotient of this fibre product by its
characteristic foliation. In general, even if the fibre product is
smooth, $S*S'$ is still not a smooth manifold, but just a quotient
of a manifold by a foliation.  We will have to deal with this
problem later, see Remark \ref{rem:smorita}.  But when the
characteristic foliation is simple, $S*S'$ is a symplectic
manifold. We may write $S*_PS'$ instead of $S*S'$ to identify the
Poisson manifold over which we are taking the tensor product.

If one is given two left modules (one could do the same for right
modules, of course), one can apply the tensor product construction by
changing the ``handedness'' of one of them.  Thus, if $S$ and $S'$ are
left $P$-modules, then $\overline{S'}$ is a right module, and
 we can form the tensor product $\overline{S'}*{S}$. We call this the
 \textbf{classical intertwiner space} \cite{Xu91b,Xu94} of $S$ and $S'$ and denote it by
 $\mathrm{Hom}(S,S')$.  The name and notation come from the case of
 modules over
 an algebra, where the tensor product  $\Y ^*\otimes \X$
is naturally isomorphic to the space of module homomorphisms from
$\Y$ to $\X$ when these modules are ``finite dimensional''.  When
the algebra is a group algebra, the modules are representations of
the group, and the module homomorphisms are known as intertwining
operators.

\begin{example}(\textit{Symplectic reduction})\label{ex:tensor-red}

Let $J:S \to \mathfrak{g}^*$ be the momentum map for a hamiltonian
action of a connected Lie group $G$ on a symplectic manifold $S$.  Let
$S'=\mathcal{O}_{\mu}$ be the coadjoint orbit through $\mu \in
\mathfrak{g}^*$, equipped with the symplectic structure induced by the
Lie-Poisson structure on $\mathfrak{g}^*$, and let $\iota:
\mathcal{O}_{\mu} \hookrightarrow \mathfrak{g}^*$ be the inclusion,
which is a Poisson map.  Then the classical intertwiner space
$\homom(S,\calo_\mu)$ is equal to
$J^{-1}(\mathcal{O}_{\mu})/G \cong J^{-1}(\mu)/G_\mu,$
i.e., the symplectic reduction of $S$ at the momentum value $\mu$.
\end{example}

A ($P_1$,$P_2$)-\textbf{bimodule} is a symplectic
manifold $S$ and a pair of
maps $\PSP$ making $S$ into a left $P_1$-module and a right
$P_2$-module and satisfying the ``commuting actions'' condition:
\begin{equation}\label{eq:commact}
\{J_1^*C^\infty(P_1),J_2^*C^\infty(P_2)\}=0.
\end{equation}
(Such geometric bimodules, without the completeness assumption,
 are called \textbf{dual pairs} in \cite{We83}.)
An \textbf{isomorphism} of bimodules is a symplectomorphism
commuting with the Poisson maps.

Given bimodules $\PSP$ and $\PSPprime$, we may form the tensor product
$S*_{P_2}S'$, and it is easily seen that this tensor product, whenever it is smooth,
becomes a ($P_1,P_3$)-bimodule \cite{Xu91b,Land01}.  We think of this tensor product as the
\textbf{composition} of $S$ and $S'$.

\begin{remark}(\textit{Modules as bimodules and geometric Rieffel Induction})\label{rmk-modbimod}

For any left $P_2$-module $S'$, there is an associated bimodule
$P_2\leftarrow S' \rightarrow \pt$, where $\pt$ is just a point.  Given
a bimodule $\PSP$, we can form its tensor product with $P_2\leftarrow
S' \rightarrow \pt$ to get a $(P_1,\pt)$-bimodule.  In this way, the
$(P_1,P_2)$-bimodule ``acts'' on $P_2$-modules to give $P_1$-modules.
This is the geometric analogue of the functors \eqref{eq:alginducton} and
\eqref{eq:rieffelinduction},  for unital and $C^*$-algebras, respectively.
\end{remark}

\begin{example}\label{ex:tensor-red2}
Following Example \ref{ex:tensor-red}, suppose that the orbit space
$S/G$ is smooth, in which case it is a Poisson manifold in a natural
way. Consider the bimodules $S/G \leftarrow S
\stackrel{J}{\to}
\overline{\mathfrak{g}^*}$ and $\overline{\mathfrak{g}^*}
\stackrel{\iota}{\leftarrow}
\overline{\mathcal{O}_{\mu}} \to \mathrm{pt}$.
Their tensor product is the $(S/G,\mathrm{pt})$-bimodule
$S/G \leftarrow S*\overline{\mathcal{O}_{\mu}}\rightarrow \mathrm{pt}$,
where the map on the left is the inclusion of the symplectic
reduced space $S*\overline{{\calo}_{\mu}}=\homom(S,\calo)$ as a symplectic leaf of
$S/G$.
\end{example}

Following the analogy with algebras, it is natural to think of
isomorphism classes of bimodules as generalized morphisms of
Poisson manifolds. The extra technical difficulty in this
geometric context is that tensor products do not always result in
smooth spaces. So one needs a suitable notion of ``regular
bimodules'', satisfying extra regularity conditions to guarantee
that their tensor products are smooth and again ``regular'', see
\cite{BuWe,Land01}, or an appropriate notion of bimodule modeled
on ``singular'' spaces. We will come back to these topics in
Section \ref{subsec:moritainv}.

\subsection{Symplectic groupoids}\label{subsec:sg}

In order to regard geometric bimodules over Poisson manifolds as
morphisms in a category, one needs to identify the bimodules which
serve as identities, i.e., those satisfying
$$
S*S'\cong S' \;\; \mbox{ and } \;\; S''*S \cong S''
$$ for any other bimodules $S'$ and $S''$. As we saw in Section
\ref{subsec:rings}, in the case of unital algebras, the identity
bimodule of an object $\A$ in $\catalg$ is just $\A$ itself, regarded
as an $(\A,\A)$-bimodule in the usual way.  This idea cannot work for
Poisson manifolds, since they are generally not symplectic, and
because we do not have commuting left and right actions of $P$ on
itself.
Instead,
it is the symplectic groupoids \cite{We87} which serve as such
``identity bimodules'' for Poisson manifolds, see \cite{Land01}.
If $\PGPst$ is an identity bimodule for a Poisson manifold $P$, then
there exists, in particular, a symplectomorphism $\grd * \grd \to
\grd$, and the composition
$$
\grd\times_{(s,t)} \grd \to \grd*\grd \stackrel{\sim}{\to} \grd
$$
defines a map $m: \grd\times_{(s,t)} \grd \to \grd$ which turns
out to be a groupoid multiplication\footnote{For expositions on
groupoids, we refer to \cite{SilWein99,MM03,MM04}; we adopt the
convention that, on a Lie groupoid $\grd$ over $P$, with source
$s$ and target $t$, the multiplication is defined on $\{(g,h)\in
\grd\times\grd,\, s(g)=t(h)\}$, and we identify the Lie algebroid
$A(\grd)$ with $\ker(Ts)|_P$, and $Tt$ is the anchor map.  The bracket
on the Lie algebroid comes from identification with
\textit{right}-invariant vector fields, which is counter to a
convention often used for Lie groups.},
compatible with the symplectic form on $\grd$ in the sense that
$\gra(m) \subseteq \grd \times \overline{\grd} \times
\overline{\grd}$ is a \textit{lagrangian} submanifold. If $p_i:
\grd\times_{(s,t)}\grd \to \grd$, $i=1,2$, are the natural
projections, then the compatibility between $m$ and $\omega$ is
equivalent to the condition
\begin{equation}\label{eq:multipl}
m^*\omega= p_1^*\omega + p_2^*\omega.
\end{equation}
A 2-form $\omega$ satisfying \eqref{eq:multipl} is called \textbf{multiplicative} (note that if
$\omega$ were a function, \eqref{eq:multipl} would mean that
$\omega(gh)=\omega(g) + \omega(h)$), and a groupoid equipped with a
multiplicative symplectic form
is called a \textbf{symplectic groupoid}.

If $(\grd,\omega)$ is a symplectic groupoid over a manifold $P$, then the
following important properties follow from
the compatibility condition \eqref{eq:multipl}, see \cite{CDW87}:

\begin{enumerate}
\item[$i)$] The unit section $P \hookrightarrow \grd$ is lagrangian;
\item[$ii)$] The inversion map $\grd \to \grd$ is an anti-symplectic involution;
\item[$iii)$] The fibres of the target and source maps, $t, s:\grd \to P$, are the symplectic orthogonal
of one another;
\item[$iv)$] At each point of $\grd$, $\ker(Ts)=\{X_{t^*f}\,|\, f\in C^\infty(P)\}$ and
$\ker(Tt)=\{X_{s^*f}\,|\, f\in C^\infty(P)\}$;
\item[$v)$] $P$ carries a unique Poisson structure such that the target map $t$ is a
Poisson map (and the source map $s$ is anti-Poisson).
\end{enumerate}
A Poisson manifold $(P,\Pi)$ is called \textbf{integrable} if there exists a symplectic groupoid
$(\grd,\omega)$ over $P$ which induces $\Pi$ in the sense of $v)$, and we refer to $\grd$ as an
\textbf{integration} of $P$. As we will discuss later, not every Poisson manifold
is integrable is this sense, see \cite{CrFe02,We87}.
But if $P$ is integrable, then there exists a symplectic groupoid
integrating it which has simply-connected (i.e., connected with
trivial fundamental group) source fibres \cite{MaXu}, and this groupoid is unique up
to isomorphism.

\begin{remark}(\textit{Integrability and complete symplectic realizations})\label{rem:complete}

If $(\grd,\omega)$ is an integration of $(P,\Pi)$, then the target
map $t:\grd \to P$ is a Poisson submersion which is always
complete. On the other hand, as proven in \cite{CrFe02}, if a
Poisson manifold $P$ admits a complete symplectic realization $S
\to P$ which is a submersion, then $P$ must be integrable.
\end{remark}

\begin{remark}(\textit{The Lie algebroid of a Poisson manifold})

All the integrations of a Poisson manifold $(P,\Pi)$ have (up to
natural isomorphism) the same Lie algebroid.  It is $T^*P$, with a
Lie algebroid structure with anchor $\widetilde{\Pi}:T^*P\to TP$,
and Lie bracket on $\Gamma(T^*P)=\Omega^1(P)$ defined by
\begin{equation}\label{eq:brk1form}
[\alpha,\beta]:= \Lie_{\widetilde{\Pi}(\alpha)}(\beta) -
\Lie_{\widetilde{\Pi}(\beta)}(\alpha) -d\Pi(\alpha,\beta).
\end{equation}
Note that \eqref{eq:brk1form} is uniquely characterized by $[df,dg]=d\{f,g\}$ and the Leibniz identity.
Following Remark \ref{rem:algebroids}, we know that  $L_{\Pi}=\gra(\widetilde{\Pi})$
also carries a Lie algebroid structure, induced by the Courant bracket.
The  natural projection
$\pr_2:TP\oplus T^*P \to T^*P$ restricts to a vector bundle isomorphism $L_{\Pi} \to T^*P$
which defines an isomorphism of Lie algebroids.

On the other hand, if $(\grd,\omega)$ is a symplectic groupoid integrating $(P,\Pi)$,
then the bundle isomorphism
\begin{equation}\label{eq:bundmap}
\Ker(Ts)|_P \longrightarrow T^*P,\;\; \xi \mapsto i_{\xi}\omega|_{TP}
\end{equation}
induces an isomorphism of Lie algebroids $A(\grd) \stackrel{\sim}{\to} T^*P$, where
$A(\grd)$ is the Lie algebroid of $\grd$, so the symplectic groupoid
$\grd$ integrates $T^*P$ in the sense of Lie algebroids.
It follows from \eqref{eq:bundmap} that $\dim(\grd)=2\dim(P)$.
\end{remark}

In the work of Cattaneo and Felder \cite{CaFe}, symplectic
groupoids arise as reduced phase spaces of Poisson
sigma models. This means that one begins with the space of
paths on $T^*P$, which has a natural symplectic structure,
restricts to a certain submanifold of ``admissible'' paths, and
forms the symplectic groupoid $\grd(P)$ as a quotient of this submanifold
by a foliation. This can also be described as an
infinite-dimensional symplectic reduction. The resulting space
is a groupoid but may not be a manifold.  When it is a manifold, it is
a the source-simply-connected symplectic groupoid of $P$.
When
$\grd(P)$ is not a manifold, as the leaf space of a foliation, it can
be considered as a differentiable stack, and even as a symplectic
stack.   In the world of stacks \cite{Metzler}, it is again a smooth groupoid; we
will call it an \textbf{S-groupoid}.  The first steps of this program have been
carried out by Tseng and Zhu  \cite{TZ}.  (See \cite{We04} for an
exposition, as well as Remark \ref{rem:smorita} below.)

This construction of symplectic groupoids has been extended to
general Lie algebroids, see \cite{CrFe01,Se}.
Crainic and Fernandes  \cite{CrFe01} describe explicitly the
obstructions to the integrability of Lie algebroids and, in
\cite{CrFe02}, identify these obstructions for the case of Poisson
manifolds and symplectic groupoids.  Integration by S-groupoids is
done in \cite{TZ}.

The next three examples illustrate simple yet important
 classes of integrable
Poisson manifolds and their symplectic groupoids.

\begin{example}(\textit{Symplectic manifolds})

If $(P,\omega)$ is a symplectic manifold, then the pair groupoid
$P \times P$ equipped with the symplectic form $\omega \times
(-\omega)$ is a symplectic groupoid integrating $P$. In order to
obtain a source-simply-connected integration, one should consider
the fundamental groupoid $\pi(P)$, with symplectic
structure given by the pull-back of the symplectic form on
$P\times \overline{P}$ by the covering map $\pi(P) \to
P\times \overline{P}$.
\end{example}

\begin{example}(\textit{Zero Poisson brackets})\label{ex:zero}

If $(P,\Pi)$ is a Poisson manifold with $\Pi=0$, then
$\grd(P)=T^*P$. In this case, the source and target maps coincide
with the projection $T^*P\to P$, and the multiplication on $T^*P$
is given by fibrewise addition. There are, however, other
symplectic groupoids integrating $P$, which may not have connected
or simply-connected source fibres.  For example, if $T^*P$ admits
a basis of closed $1$-forms, we may divide the fibres of $T^*P$ by the
lattice generated by these forms to obtain a groupoid whose source and
target fibres are tori.  Or, if $P$ is just a point, any discrete
group is a symplectic groupoid for $P$.  We refer to \cite{BuWe} for
more details.
\end{example}

\begin{example}(\textit{Lie-Poisson structures})\label{ex:lie-poiss-grd}

Let $P=\mathfrak{g}^*$ be the dual of a Lie algebra $\mathfrak{g}$,
equipped with its Lie-Poisson structure, and let $G$ be a Lie group
with Lie algebra $\mathfrak{g}$. The transformation groupoid
$G \ltimes \mathfrak{g}^*$ with respect to the coadjoint action,
equipped with the symplectic form obtained from the identification
$G \times \mathfrak{g}^* \cong T^*G$ by right translation, is
a symplectic groupoid integrating $\mathfrak{g}^*$.
This symplectic groupoid is  source-simply-connected
just when $G$ is a (connected)  simply-connected Lie group.
\end{example}

\begin{remark}(\textit{Lie's third theorem})

Let $\mathfrak{g}$ be a Lie algebra.
Example \ref{ex:lie-poiss-grd} shows that integrating
$\mathfrak{g}$, in the usual sense of finding a Lie group $G$ with Lie algebra
$\mathfrak{g}$, yields an integration of the Lie-Poisson structure of $\mathfrak{g}^*$.
On the other hand, one can use the integration of the
Lie-Poisson structure of $\mathfrak{g}^*$ to construct a
Lie group integrating $\mathfrak{g}$. Indeed, if $\grd$ is a
symplectic groupoid integrating $\mathfrak{g}^*$, then the map
$$
\mathfrak{g} \to \mathcal{X}(\grd),\;\; v \mapsto X_{t^*v}
$$
is a faithful representation of $\mathfrak{g}$ by vector fields on
$\grd$. Here $t:\grd \to \mathfrak{g}^*$ is the target map, and we
regard $v \in \mathfrak{g}$ as a linear function on
$\mathfrak{g}^*$. We then use the flows of these vector fields to
define a (local) Lie group integrating $\mathfrak{g}$.  If
we fix $x \in \grd$, the ``identity'' of the local Lie group, so
that $t(x)=0$, then the Lie group sits in $\grd$ as a lagrangian
subgroupoid.  So the two ``integrations'' are the same.

The idea of using a symplectic realization of $\mathfrak{g}^*$ to find
a Lie group integrating $\mathfrak{g}$ goes back to Lie's original proof of
``Lie's third theorem.''
A regular point of $\mathfrak{g}^*$ has
a neighborhood $U$ with
coordinates $(q_1,\ldots,q_k,p_1,\ldots,p_k, e_1,\ldots,e_l)$ such that the
Lie-Poisson structure can be written as
$$
\sum_{i=1}^k \frac{\partial}{\partial{q_i}}\wedge \frac{\partial}{\partial{p_i}}
$$
(see Section \ref{subsec:local}). The map  $\mathfrak{g} \to \mathcal{X}(U)$,
$v \mapsto -X_v$ is a Lie algebra homomorphism, but not faithful in general.
It suffices, though, to
add $l$ new coordinates $(f_1,\ldots,f_l)$ and consider the local
symplectic realization $U \times \mathbb{R}^l \to U$, with symplectic Poisson structure
$$
\Pi'=\sum_{i=1}^k \frac{\partial}{\partial{q_i}}\wedge \frac{\partial}{\partial{p_i}}
+\sum_{i=1}^l \frac{\partial}{\partial{e_i}}\wedge \frac{\partial}{\partial{f_i}}.
$$
The map $\mathfrak{g} \to \mathcal{X}(U\times \mathbb{R}^l)$, $v
\mapsto -X'_v:=\widetilde{\Pi}'(v)$, is now a faithful Lie algebra
homomorphism.  Once again, we can use the flows of the hamiltonian
vector fields of the coordinates on $\frakg$ to construct a local
Lie group.

More generally, if $\grd$ is a Lie groupoid and $A$ is its Lie algebroid,
then $T^*\grd$ is naturally a symplectic groupoid over $A^*$, see \cite{CDW87}. The
induced Poisson structure on $A^*$ is a generalization of a
Lie-Poisson structure. Conversely, if $A$ is an integrable Lie
algebroid, then $\grd(A)$, its source-simply-connected integration,
can be constructed as a lagrangian subgroupoid of the symplectic groupoid $\grd(A^*)$ integrating
$A^*$ \cite{Catt}.
\end{remark}

The following is an example of a non-integrable Poisson structure.

\begin{example}\label{ex:noninteg}(\textit{Nonintegrable Poisson structure})

Let $P=S^2 \times \mathbb{R}$. Let $\Pi_{S^2}$ be the natural
symplectic structure on $S^2$. Then the product Poisson structure
on $P$, $\Pi_{S^2}\times \{0\}$ is integrable. But if we multiply
this Poisson structure by $(1+t^2)$, $t \in \mathbb{R}$ (or use
any other nonconstant function which has a critical point), then
the resulting Poisson structure $(1+t^2)(\Pi_{S^2}\times \{0\})$
is \textit{not} integrable \cite{CrFe02,We87}. In this case the
symplectic S-groupoid $\grd(P)$ is not a manifold.

We will have more to say about this example in Section
\ref{subsec:moritainv}.
\end{example}

\begin{remark}(\textit{Twisted presymplectic groupoids})

Let $\grd$ be a Lie groupoid over a manifold $P$. For each $k > 0$, let $\grd_k$
be the manifold of composable sequences of $k$-arrows,
$$
\grd_k := \grd \times_{(s,t)} \grd \times_{(s,t)} \dots
\times_{(s,t)}\grd,\qquad \mbox{ ($k$ times)}
$$
and set $\grd_0 = P$. The sequence of manifolds $\grd_k$, together with the
natural maps  $\partial_i: \grd_{k} \to \grd_{k-1}$, $i=0,\ldots,k$,
$$
\partial_i(g_1,\dots,g_k)=
\left \{
\begin{matrix}
(g_2,\ldots,g_k), & \mbox{if} & i=0,\\
(g_1,\ldots,g_ig_{i+1},\ldots,g_k), &\mbox{if} & 0<i<k\\
(g_1,\ldots,g_{k-1}) &\mbox{if} & i=k.
\end{matrix}
\right .
$$
defines a simplicial manifold $\grd_{\bullet}$. The \textbf{bar-de Rham complex} of $\grd$ is
the total complex of  the double
complex $\Omega^\bullet(\grd_{\bullet})$, where the boundary maps are $
d:\Omega^{q}(\grd_k) \to \Omega^{q+1}(\grd_k),
$
the usual de Rham differential, and
$
\partial : \Omega^q(\grd_k) \to \Omega^{q}(\grd_{k+1}),
$
the alternating sum of the pull-back of the $k+1$ maps $\grd_k \to \grd_{k+1}$, as in group cohomology.
For example, if $\omega \in \Omega^2(\grd)$, then
$$
\partial \omega = p_1^*\omega -m^*\omega + p_2^*\omega.
$$
(As before, $m$ is the groupoid multiplication, and $p_i:\grd_2
\to \grd$, $i=1,2$, are the natural projections.) It follows that
a  2-form $\omega$ is a 3-cocycle in the total complex if and only
if it is  multiplicative and closed; in particular, a symplectic
groupoid can be defined as  a Lie groupoid $\grd$ together with a
nondegenerate 2-form $\omega$ which is a 3-cocycle.

More generally, one can consider 3-cochains which are  sums
$\omega+\phi$, where $\omega \in \Omega^2(\grd)$ and $\phi \in \Omega^3(P)$. In this case,
the coboundary condition is that $d\phi=0$, $\omega$ is multiplicative, and
$$
d\omega = s^*\phi - t^*\phi.
$$ A groupoid $\grd$ together with a 3-cocycle $(\omega,\phi)$ such
that $\omega$ is nondegenerate is called a \textbf{$\phi$-twisted
symplectic groupoid} \cite{SeWe01}.  Just as symplectic groupoids are
the global objects associated with Poisson manifolds, the
twisted symplectic groupoids are the global objects associated with twisted
Poisson manifolds \cite{CaXu}.

Without non-degeneracy assumptions on $\omega$, one has the following result concerning
the infinitesimal version of 3-cocycles \cite{BCWZ}:
If $\grd$ is source-simply connected and $\phi\in \Omega^3(P)$, $d\phi=0$, then there is
a one-to-one correspondence between 3-cocycles $\omega + \phi$ and bundle maps
$\sigma: A \to T^*P$ satisfying the following two conditions:
\begin{eqnarray}
\SP{\sigma(\xi),\rho(\xi')}&=&-\SP{\sigma(\xi'),\rho(\xi)};\\\label{eq:cond1}
\sigma([\xi,\xi'])&=& \Lie_\xi(\sigma(\xi'))-\Lie_{\xi'}(\sigma(\xi))+ d\SP{\sigma(\xi),\rho(\xi')}+
i_{\rho(\xi)\wedge\rho(\xi'))}(\phi),\label{eq:cond2}
\end{eqnarray}
where $A$ is the Lie algebroid of $\grd$, $[\cdot,\cdot]$ is the bracket on $\Gamma(A)$,
$\rho:A \to TP$ is the anchor, and $\xi,\xi' \in \Gamma(A)$. For one direction of this
correspondence, given
$\omega$, the associated bundle map $\sigma_{\omega}:A \to T^*P$ is just $\sigma_{\omega}(\xi)=
i_\xi\omega|_P$.

For a given $\sigma:A\to T^*P$ satisfying \eqref{eq:cond1}, \eqref{eq:cond2},
let us consider the bundle map
\begin{equation}\label{eq:rhosigma}
(\rho,\sigma): A \to TP\oplus T^*P.
\end{equation}
A direct computation shows that if the rank of $L_{\sigma}:=\mathrm{Image}(\rho,\sigma)$
equals $\mathrm{dim}(P)$, then $L_{\sigma} \subset TP\oplus T^*P$  is a $\phi$-twisted
Dirac structure on $P$. In this case, it is easy to check that  \eqref{eq:rhosigma}
yields a (Lie algebroid) isomorphism $A \to L_{\sigma}$ if and only if
\begin{itemize}
\item[1)] $\dim(\grd)=2\dim(P)$;
\item[2)] $\ker(\omega_x)\cap\ker(T_xs)\cap\ker(T_xt)=\{0\}$\; \mbox{for all $x \in P$.}
\end{itemize}
A groupoid $\grd$ over $P$ satisfying $1)$ together with a 3-cocycle $\omega+\phi$ so that $\omega$
satisfies $2)$ is called a \textbf{$\phi$-twisted presymplectic groupoid} \cite{BCWZ,Xu03}.
As indicated by the previous discussion, they are precisely the global objects
integrating twisted Dirac structures. The 2-form $\omega$ is nondegenerate
if and only if the associated Dirac structure is Poisson, recovering the known
correspondence between (twisted) Poisson structures and (twisted) symplectic groupoids.
\end{remark}

The following example describes presymplectic groupoids integrating Cartan-Dirac
structures; it is analogous to Example \ref{ex:lie-poiss-grd}.

\begin{example} {\it(Cartan-Dirac structures and the
AMM-groupoid)}\label{ex:AMMgrd}

Let $G$ be a Lie group with Lie algebra $\mathfrak{g}$, equipped with a
nondegenerate bi-invariant quadratic form $(\cdot,\cdot)_{\mathfrak{g}}$.
The {\bf AMM groupoid} \cite{BeXuZh} is the action groupoid $\grd =G \ltimes G$ with respect to the
conjugation action, together with the 2-form \cite{AMM}
$$
\omega_{(g, x)} = \frac{1}{2} \left((\mathrm{Ad}_{x} p_g^* \lambda, p_g^*
\lambda)_{\mathfrak{g}} +
( p_g^*\lambda, p_x^*(\lambda+\bar{\lambda}))_{\mathfrak{g}}\right),
$$
where $p_g$ and $p_x$ denote the projections onto the
first and second
components of $G\times G$, and $\lcart$ and $\rcart$ are the left and right Maurer-Cartan
forms. The AMM-groupoid is a $\phi^G$-twisted presymplectic groupoid integrating
$L_G$ \cite{BCWZ}, the Cartan-Dirac
structure on $G$ defined in  Example \ref{ex:cartan-dirac}.
If $G$ is simply connected, then  $(G\ltimes G, \omega)$ is isomorphic to $\grd(L_G)$, the
source-simply connected integration of $L_G$;
in general, one must pull-back $\omega$ to $\widetilde{G}\ltimes G$, where
$\widetilde{G}$ is the universal cover of $G$.
\end{example}


\subsection{Morita equivalence for groups and groupoids}\label{subsec:moritagr}
Since groupoids play such an important role in the Morita equivalence of
Poisson manifolds, we will take some time to discuss Morita
equivalence of groupoids in general.  We begin with groups.

If we try to define Morita equivalence of groups as equivalence
between their (complex linear) representation categories, then we are
back to algebra, since representations of a group are the same as
modules over its group algebra over $\mathbb{C}$.  (This is
straightforward for discrete groups, and more elaborate for
topological groups.)  Here, we just remark that nonisomorphic groups
can have isomorphic group algebras (e.g. two finite abelian groups
with the same number of elements), or more generally Morita equivalent
group algebras (e.g. two finite groups with the same number of
conjugacy classes, hence the same number of isomorphism classes of
irreducible representations).

We obtain a more geometric notion of Morita equivalence for groups by
considering actions on manifolds rather than on linear spaces.  Thus, for
Lie groups (including discrete groups) $G$ and $H$, bimodules
are $(G,H)$-``bispaces'', i.e.
manifolds where $G$ acts on the left, $H$ acts on the right, and
the actions commute. The ``tensor product'' of such bimodules is
defined by the orbit space
$$
\GXH * \HYK := \X \times \Y/H,
$$
where $H$ acts on $\X \times \Y$ by $(x,y) \mapsto (xh,h^{-1}y)$.
The result of this
operation may no longer be smooth, even if $\X$ and $\Y$ are.
Under suitable regularity assumptions, to be explained below,
the tensor product is a smooth manifold, so we
consider the category in
which objects are groups and morphisms are isomorphism classes of ``regular'' bispaces, and
we define
\textbf{Morita equivalence} of groups as isomorphism in this category. Analogously
to the case of algebras, we have an associated notion of \textbf{Picard group(oid)}.

\begin{exercise}
Show that a bispace $\GXH$ is ``invertible'' if and
only if the $G$ and $H$-actions are free and transitive.
\end{exercise}

If $\GXH$ is invertible and
we fix a point $x_0 \in \X$, by the result of the previous
exercise,
there exists for each $g \in G$
a unique $h \in H$ such that $gx_0h^{-1}=x_0$. The correspondence
$g \mapsto h$ in fact establishes a group isomorphism $G \to H$.
So, for groups, Morita equivalence induces the same equivalence
relation as the usual notion of isomorphism.  As we will see in Example \ref{ex:picgrgr}
of Section \ref{subsec:picard}, the situation for Picard groups resembles
somewhat that for algebras, where outer automorphisms play a key role.

For a full discussion of Morita equivalence of Lie groupoids, we refer
to the article of Moerdijk and Mr\v{c}un \cite{MM04} in this volume.
Here, we will briefly summarize the theory.

An action (from the left) of a Lie groupoid $\grd$ over $P$ on a
manifold $S$ consists of a map $J:S\to P$
and a map $\grd \times_{(s,J)}S \to S$ (where $s$ is the source
map of $\grd$) satisfying axioms analogous to those of a group
action;  $J$ is sometimes  called the \textbf{moment}
of the action (see Example \ref{ex:hamilt}).
The action is  \textbf{principal} with respect to a map $p:S \to
M$ if $p$ is a surjective submersion and if $\grd$ acts freely and
transitively on each $p$-fibre; principal $\grd$-bundles
are also called \textbf{$\grd$-torsors}.

Right actions and torsors are defined in the obvious analogous way.
If groupoids $\grd_1$ and $\grd_2$
act on $S$ from the left and right, respectively, and the actions
commute, then we call $S$
a \textbf{$(\grd_1,\grd_2)$-bibundle}. A bibundle is \textbf{left
principal} when the left $\grd_1$-action is principal with respect
to the moment map for the right action of $\grd_2$.

If $S$ is a $(\grd_1,\grd_2)$-bibundle with moments $\PSP$, and if
$S'$ is a $(\grd_2,\grd_3)$-bibundle with moments
$P_2\stackrel{J_2'}{\leftarrow}S'\stackrel{J_3'}{\to}P_3$, then
their ``tensor product'' is the orbit space
\begin{equation}\label{eq:grdtensor}
S*S':=(S\times_{(J_2,J_2')}S')/\grd_2,
\end{equation}
where $\grd_2$ acts on $S\times_{(J_2,J_2')}S'$ diagonally. The assumption that $S$ and
$S'$ are left principal guarantees that $S*S'$ is a smooth
manifold and that its natural $(\grd_1,\grd_3)$-bibundle structure
is left principal.

Two $(\grd_1,\grd_2)$-bibundles are \textbf{isomorphic} if there
is a diffeomorphism between them commuting with the groupoid
actions and their moments. The ``tensor product'' \eqref{eq:grdtensor} is
associative up to natural isomorphism, so we may define a category
$\catgr$ in which the objects are Lie groupoids and morphisms are
isomorphism classes of left principal bibundles. Just as in the
case of algebras, we call two Lie groupoids \textbf{Morita
equivalent} if they are isomorphic as objects in $\catgr$, and we
define the associated notion of Picard group(oid) just as we do
for algebras. We note that a $(\grd_1,\grd_2)$-bibundle $S$ is
``invertible'' in $\catgr$ if and only if it is
\textbf{biprincipal}, i.e., principal with respect to both left
and right actions; a biprincipal bibundle is also called a
\textbf{Morita equivalence} or a \textbf{Morita bibundle}.

\begin{example}(\textit{Transitive Lie groupoids})\label{ex:transitive}

Let $\grd$ be a Lie groupoid over $P$.
For a fixed $x\in P$, let $\grd_x$ be the isotropy group of $\grd$ at $x$,
and let $E_x=s^{-1}(x)$. Then $E_x$ is a $(\grd,\grd_x)$-bibundle.
It is a Morita bibundle if and only if $\grd$ is transitive, i.e.,
for any $x, y \in P$, there exists $g\in \grd$ so that $s(g)=y$ and $t(g)=x$.
In fact, a Lie groupoid is transitive if and only if it is Morita equivalent to a Lie group.
\end{example}

\subsection{Modules over Poisson manifolds and symplectic groupoid actions}\label{subsec:actions}

Example \ref{ex:gstarmod} shows that modules over $\mathfrak{g}^*$ are the same thing
as hamiltonian $G$-manifolds, where $G$ is the connected and simply connected Lie group
with Lie algebra $\mathfrak{g}$. As we discuss in this section, this is a
particular case of a much more general correspondence between
modules over Poisson manifolds and symplectic groupoid actions.

Let $(\grd,\omega)$ be a symplectic groupoid over $P$ acting
on a symplectic manifold $(S,\omega_S)$ with moment $J$.
Let $a: \grd \times_{(s,J)} S \to S$ denote the action.
We call the action \textbf{symplectic} if it satisfies
the property (analogous to the condition on multiplicative forms)
that $\gra(a) \subset \grd(P) \times S \times \overline{S}$ is
lagrangian. Equivalently, $a$ is symplectic if
\begin{equation}\label{eq:compactact}
a^*\omega_S = p_{\st{S}}^*\omega_S + p_{\st{\grd}}^*\omega,
\end{equation}
where $p_{\st{S}}:\grd \times_{(s,J)} S \to S$ and $p_{\st{\grd}}:
\grd \times_{(s,J)} S \to \grd$ are the natural projections.

A key observation relating actions of symplectic groupoids to modules over Poisson manifolds
is that if $J:S\to P$ is the moment map for a \textit{symplectic} action of
a symplectic groupoid $\grd$ over $P$,
then $J$ is automatically a complete Poisson
map \cite{MiWe}, defining a module over $P$. On the other hand,
a module $J:S\to P$ over an
integrable Poisson manifold $P$ automatically carries a \textit{symplectic} action of
the source-simply connected symplectic groupoid $\grd(P)$. So there is a
one-to-one correspondence between $P$-modules and symplectic
actions of $\grd(P)$ \cite{CDW87}.

\begin{example}\label{ex:hamilt} (\textit{Hamiltonian spaces})

Let $G$ be a simply-connected Lie group with Lie algebra $\mathfrak{g}$. Any
complete symplectic realization $J:S \to \mathfrak{g}^*$ induces
an action of the symplectic groupoid $T^*G$ on $S$:
$$
 \xymatrix {
{T^*G} \ar@<-0.5ex>[d]\ar@<0.5ex>[d]  & S\ar[dl]^{J}\\
\mathfrak{g}^*
 }
$$
In this case, $T^*G = G \ltimes \mathfrak{g}^*$ is a
transformation Lie groupoid, and, as such, its action is
equivalent to an ordinary $G$-action on $S$ for which $J$ is
$G$-equivariant. Moreover, the $G$-action corresponding to the
symplectic $T^*G$-action induced by $J:S\to P$ is a hamiltonian
$G$-action for which $J$ is a momentum map. So we recover the
result of Example \ref{ex:gstarmod} on the isomorphism
(not only equivalence) between the categories of complete
symplectic realizations of $\mathfrak{g}^*$ and hamiltonian
$G$-manifolds.  Notice that the momentum map for the group action is
the moment map for the groupoid action; it is this example which
motivates the term ``moment'' as applied to groupoid actions.
\end{example}

\begin{remark}(\textit{Infinitesimal actions})

The relationship between complete symplectic realizations and
symplectic groupoid actions has an  infinitesimal
counterpart. A symplectic realization (not necessarily complete) $J:S
\to P$ induces a Lie algebra homomorphism
\begin{equation}\label{eq:infact}
\Omega^1(P) \to \mathcal{X}(S),\;\; \alpha \mapsto
\widetilde{\Pi}_S(J^*\alpha),
\end{equation}
where the bracket on 1-forms is the one of \eqref{eq:brk1form}.
This maps defines a Lie algebroid action
of the Lie algebroid of
$P$, $T^*P$, on $S$.
The completeness of $J$ allows this infinitesimal action to be
integrated to an action of the source-simply-connected integration
$\grd(P)$ (see \cite{MM02}),
and this action turns out to be symplectic.
\end{remark}

\begin{remark}(\textit{Symplectic groupoid actions on Poisson
    manifolds})\label{rem:hamiltpoiss}

As in Remark \ref{rem:genmod}, the correspondence between
$P$-modules and symplectic $\grd(P)$-actions holds in more
generality: a Poisson map $Q \to P$ from \textit{any} Poisson
manifold $Q$ induces an infinitesimal $T^*P$-action on $Q$, by the
same formula as in \eqref{eq:infact}. When the Poisson map is
complete (and $P$ is integrable), it gives rise to an action of
$\grd(P)$ on $Q$, which preserves the symplectic leaves of $Q$;
its restriction to each leaf is a symplectic action.  The action
is a Poisson action in the sense that its graph is lagrangian
\cite{We88} in the appropriate product, see \cite{BuCr} for
details.
\end{remark}

\begin{remark}(\textit{Realizations of Dirac structures and presymplectic groupoid
actions})\label{rem:diracreal}

The correspondence between modules over a Poisson manifold $P$ and
symplectic actions of $\grd(P)$ extends to one between ``modules''
over Dirac manifolds and suitable actions of presymplectic
groupoids \cite{BuCr,BCWZ}.

In order to introduce the notion of ``realization'' of a Dirac
manifold, let us note that, if $(P,\Pi)$ is a Poisson manifold,
then the infinitesimal $T^*P$-action \eqref{eq:infact} induced by
a Poisson map $J:Q\to P$ can be equivalently expressed in terms of
$L_{\Pi}$ by the Lie algebra homomorphism
$$
\Gamma(L_{\Pi})\to \mathcal{X}(Q),\;\; (X,\alpha) \mapsto Y,
$$
where $Y$ is uniquely determined by the condition
$(Y,J^*\alpha)\in L_{\Pi_Q}$. Since $J$ is a Poisson map, it also
follows that $X=TJ(Y)$.

If $(P,L)$ and $(Q,L_Q)$ are Dirac manifolds, and $J:Q\to P$ is a
forward Dirac map, then \eqref{eq:forward} implies that for each
$(X,\alpha) \in L$ over the point $J(y)\in P$, there exists $Y\in T_yQ$
such that $(Y,TJ^*(\alpha))\in (L_Q)_y$ and $X=T_yJ(Y)$. However,
unlike the situation of Poisson maps, $Y$ is \textit{not} uniquely
determined by these conditions; this is the case if and only if
\begin{equation}\label{eq:nondegcond}
\ker(TJ)\cap\ker(L_Q)=\{0\}.
\end{equation}
If \eqref{eq:nondegcond} holds at all points of $Q$, then the
induced  map $\Gamma(L)\to \mathcal{X}(Q)$, $(X,\alpha)\mapsto Y$,
defines an infinitesimal $L$-action on $Q$.

A \textbf{Dirac realization} \cite{BuCr} of a $\phi$-twisted Dirac
manifold $(P,L)$ is a forward Dirac map $J:Q \to P$, where $Q$ is
a $J^*\phi$-twisted Dirac manifold and \eqref{eq:nondegcond} is
satisfied. If $Q$ is a $J^*\phi$-twisted presymplectic manifold,
then $J$ is called a \textbf{presymplectic realization}. We call a
Dirac realization \textbf{complete} if the induced infinitesimal
action is complete (in the sense of Lie algebroid actions, see
\cite{MM02}). As in the case of Poisson maps, complete Dirac
realizations $J:Q\to P$ are the same thing as global actions of
the presymplectic groupoid $\grd(L)$ on $Q$ ``compatible'' with
$L_Q$ in a suitable way \cite{BuCr} (generalizing the conditions
in \eqref{eq:compactact} and Remark \ref{rem:hamiltpoiss}).
\end{remark}

The next example illustrates the discussion in Remark
\ref{rem:diracreal} and the connection between Dirac geometry and
group-valued momentum maps \cite{AMM,AKM}.

\begin{example}(\textit{Modules over Cartan-Dirac structures and quasi-hamiltonian actions})
\label{ex:cartdirmod}

As we saw in Example \ref{ex:gstarmod}, symplectic realizations of
(resp. Poisson maps into) the Lie-Poisson structure on
$\mathfrak{g^*}$ are the same thing as hamiltonian
$\mathfrak{g}$-actions on symplectic (resp. Poisson) manifolds; if
the maps are complete, one gets a correspondence with global
hamiltonian actions.

Analogously, let us consider a connected, simply-connected Lie
group $G$ equipped with $L_G$, the Cartan-Dirac structure
associated with a non-degenerate bi-invariant quadratic form
$(\cdot,\cdot)_{\mathfrak{g}}$. Then presymplectic realizations
into $G$ are exactly the same as quasi-hamiltonian
$\mathfrak{g}$-manifolds, and complete realizations correspond to global
quasi-hamiltonian $G$-actions (which can be seen as actions of the
AMM-groupoid of Example \ref{ex:AMMgrd}, analogously to Example
\ref{ex:hamilt}) \cite{BCWZ}. More
generally, (complete) Dirac realizations of $(G,L_G)$ correspond to
(global) hamiltonian quasi-Poison manifolds \cite{BuCr}, in analogy with Remark
\ref{rem:hamiltpoiss}.

In these examples, the realization maps are the group-valued
momentum maps.
\end{example}


\subsection{Morita equivalence of Poisson manifolds and symplectic groupoids}\label{subsec:geom-morita}

We now have all the ingredients which we need in
order to define a geometric notion of
Morita equivalence for Poisson manifolds which implies equivalence of
their module categories.

A \textbf{Morita equivalence} between Poisson manifolds $P_1$ and
$P_2$ is a $(P_1,P_2)$-bimodule $\PSP$ such that $J_1$ and $J_2$
are surjective submersions whose fibres are simply connected and
symplectic orthogonals of each other. By Remark
\ref{rem:complete}, Morita equivalence only applies to
\textit{integrable} Poisson structures. (The nonintegrable case
can be handled with the use of symplectic S-groupoids.  See Remark
\ref{rem:smorita}.) The bimodule $\PSbarP$, where $\overline{S}$
has the opposite symplectic structure, is also a Morita
equivalence, and $S$ and $\overline{S}$ satisfy
\begin{equation}\label{eq:invgeom}
S*_{\st{P_2}}\overline{S} \cong \grd(P_1),\;\;\mbox{ and } \;\;
\overline{S}*_{\st{P_1}}S \cong \grd(P_2).
\end{equation}
Since symplectic groupoids are ``identity bimodules'',
\eqref{eq:invgeom} is analogous to the invertibility of algebraic
bimodules \eqref{eq:invbim}.

Let us consider the category whose objects are complete symplectic
realizations  of an integrable Poisson manifold $P$, and morphisms
are symplectic maps between symplectic realizations commuting with
the realization maps. This category is analogous to the category
of left modules over an algebra, and we call it the
\textbf{category of modules over} $P$. If $\PSP$ is a Morita
equivalence, then the regularity conditions on the maps $J_1$ and
$J_2$ guarantee that if $S' \to P_2$ is a left $P_2$-module then
the tensor product $S*_{\st{P_2}}S'$ is smooth and defines a left
$P_1$-module \cite{Xu91b}. So one can define a functor between
categories of modules (i.e. complete symplectic realizations) just as
one does for algebras, see \eqref{eq:alginducton} and
\eqref{eq:alginducton2}, and prove that geometric Morita
equivalence implies the equivalence of
``representation'' categories \cite{Land01,Xu91b}:

\begin{theorem}\label{thm:equivcatpoiss}
If $P_1$ and $P_2$ are Morita equivalent, then they have
equivalent categories of complete symplectic realizations.
\end{theorem}

\begin{remark}(\textit{The ``category'' of complete symplectic
    realizations})\label{rem:categ}

In the spirit of the symplectic ``category'' of \cite{We82}, one
can also define a larger ``category'' of complete symplectic
realizations of $P$ by considering the morphisms between two
$P$-modules $J:S \to P$ and $J':S'\to P$  to be lagrangian
submanifolds in $S'\times_{(J',J)}\overline{S}$, see
\cite{Xu91,Xu91b}, with composition given by
composition   of
relations; the quotes in ``category'' are due to the fact that the
composition of two such morphisms yields another morphism only
under suitable transversality assumptions. Theorem
\ref{thm:equivcatpoiss} still holds in this more general setting
\cite{Xu91}. Unlike in the case of algebras, though, the converse
of Theorem \ref{thm:equivcatpoiss} does not hold in general
\cite{Xu91b}, see Remark \ref{rem:repvsmorita}. We will discuss
ways to remedy this problem in Section \ref{sec:geomrep} by
introducing yet another category of representations of $P$ (a
``symplectic category'').
\end{remark}

\begin{remark}(\textit{Classical intertwiner
spaces})\label{rem:intertwiners}

As a consequence of \eqref{eq:invgeom}, one can see that Morita equivalence,
in addition to establishing an equivalence of module categories,
preserves the classical intertwiner spaces.
\end{remark}

\begin{remark}(\textit{More general modules})\label{rem:general}

As indicated in Remarks \ref{rem:genmod} and
\ref{rem:hamiltpoiss}, from the point of view of hamiltonian
actions, it is natural to consider arbitrary complete Poisson maps
(not necessarily symplectic realizations) as modules over Poisson
manifolds. The ``action'' of $(P_1,P_2)$-bimodules on
$P_2$-modules in Remark \ref{rmk-modbimod}  naturally extends to
an action on Poisson maps $Q \to P_2$; in fact, one can think of
this more general tensor product as a leafwise version of the one
in Section \ref{subsec:geomrep}, and Theorem
\ref{thm:equivcatpoiss} still holds for these more general
``representations''. (This generalization is the analogue, in
algebra, of considering homomorphisms of an algebra into direct
sums of endomorphism algebras, rather than usual modules.)
\end{remark}

The notion of Morita equivalence of Poisson manifolds is closely
related to Morita equivalence of symplectic groupoids, which is a
refinement of the notion of Morita equivalence for Lie groupoids,
taking symplectic structures into account. If $\grd_1$ and
$\grd_2$ are symplectic groupoids, then a
$(\grd_1,\grd_2)$-bibundle is called \textbf{symplectic} if both
actions are symplectic. The ``tensor product'' of two symplectic
bibundles, as defined in Section \ref{subsec:moritagr}, is
canonically symplectic, so we may define a category $\SGrd$ in
which the objects are symplectic groupoids and morphisms are
isomorphism classes of left principal symplectic bibundles. (An
isomorphism between symplectic bibundles is required to preserve
the symplectic forms.) We call two symplectic groupoids $\grd_1$
and $\grd_2$ \textbf{Morita equivalent} \cite{Xu91b} if they are
isomorphic in $\SGrd$, i.e. if there exists a biprincipal
symplectic $(\grd_1,\grd_2)$-bibundle (see \cite{Land01}). A
\textbf{Morita equivalence} between symplectic groupoids is a
symplectic bibundle which is biprincipal.

If $P_1$ and $P_2$ are Poisson manifolds, and if $\PSP$ is a $(P_1,P_2)$-bimodule, then
we obtain a left symplectic action of the groupoid $\grd(P_1)$ and right symplectic
action of $\grd(P_2)$,
$$
\xymatrix { {\grd(P_1)} \ar@<-0.5ex>[d]\ar@<0.5ex>[d] &
S\ar[dl]^{J_1}\ar[dr]_{J_2} &
{\grd(P_2)} \ar@<-0.5ex>[d]\ar@<0.5ex>[d]\\
P_1 & & P_2
 }
$$
The property that $\{J_1^*C^\infty(P_1),J_2^*C^\infty(P_2)\}=0$
implies that the actions of $\grd(P_1)$ and $\grd(P_2)$ commute,
so that $S$ is a symplectic $(\grd(P_1),\grd(P_2))$-bibundle. We
say that a symplectic bimodule $\PSP$ is \textbf{regular} if the
associated symplectic $(\grd(P_1),\grd(P_2))$-bibundle is left
principal. The tensor product of symplectic bimodules defined in
Section \ref{subsec:geomrep} coincides with their tensor product
as symplectic bibundles. As a result, the tensor product of
regular symplectic bimodules is smooth and regular.

\begin{remark}(\textit{Regular bimodules})\label{rem:regbim}

Regular bimodules can be described with no reference to the
symplectic groupoid actions: $\PSP$ is regular if and only if
$J_1$ and $J_2$ are complete Poisson maps, $J_1$ is a submersion,
$J_2$ is a surjective submersion with simply-connected fibres, and
the $J_1$- and $J_2$-fibres are symplectic orthogonal of one
another.
\end{remark}

\begin{exercise}
Prove the equivalent formulation of regular bimodules in Remark
\ref{rem:regbim}. (Hint: this is a slight extension of
\cite[Thm.~3.2]{Xu91})
\end{exercise}

We define the category $\catpoiss$ in which the objects are
integrable Poisson manifolds and morphisms are isomorphism classes
of regular symplectic bimodules.

If $\PSP$ is a Morita equivalence of Poisson manifolds, then the
regularity assumptions on the maps $J_1$ and $J_2$ insure that $S$
is biprincipal for the induced actions of $\grd(P_1)$ and
$\grd(P_2)$, so that $S$ is also a Morita equivalence for the
symplectic groupoids $\grd(P_1)$ and $\grd(P_2)$. On the other
hand, if $\grd_1$ and $\grd_2$ are source-simply-connected
symplectic groupoids over $P_1$ and $P_2$, respectively, then a
$(\grd_1,\grd_2)$-Morita equivalence is a $(P_1,P_2)$-Morita
equivalence. So two integrable Poisson manifolds $P_1$ and $P_2$
are Morita equivalent if and only if their source-simply-connected
integrations, $\grd(P_1)$ and $\grd(P_2)$, are Morita equivalent
as symplectic groupoids.

\begin{remark}(\textit{Lie functor})

It follows from the discussion above that there exists a natural
equivalence between the category of
\textit{source-simply-connected} symplectic groupoids with
morphisms being Morita equivalences (resp. left principal symplectic
bibundles), and the category of
integrable Poisson manifolds with morphisms being Morita
equivalences (resp. regular bimodules).
These equivalences are similar to the one between the categories
of Lie algebras and simply-connected Lie groups, with their usual morphisms.
\end{remark}

\begin{example}(\textit{Symplectic manifolds})\label{ex:sym}

Let $P$ be a connected symplectic manifold. The universal cover
of $P$ with base point $x$, denoted $\widetilde{P}$, is a Morita
equivalence between the symplectic groupoid $\grd(P)$,
which in this case is  the fundamental groupoid of $P$, and $\pi_1(P,x)$:
\begin{equation}\label{eq:fundam}
 \xymatrix { {\grd(P)}
\ar@<-0.5ex>[d]\ar@<0.5ex>[d] & \widetilde{P}\ar[dl]\ar[dr] &
{\pi_1(P,x)} \ar@<-0.5ex>[d]\ar@<0.5ex>[d]\\
P & & {\{x\}}
 }
\end{equation}
Note that $\pi_1(P,x)$ is a symplectic groupoid for the
zero-dimensional Poisson manifold $\{x\}$,
though generally not the source-simply-connected one.

In analogy with Example \ref{ex:hamilt} on hamiltonian
actions: there is an equivalence of categories between complete
symplectic realizations of $P$ and symplectic actions of
$\pi_1(P,x)$. This suggests the slogan that ``a (connected)
 symplectic manifold $P$ with
fundamental group $\pi_1(P)$ is the dual of the Lie algebra of
$\pi_1(P)$''.

It follows from the Morita equivalence \eqref{eq:fundam} and the
discussion about Morita equivalence of groups in Section
\ref{subsec:rings} that connected symplectic manifolds $P_1$ and
$P_2$ are Morita equivalent if and only if $\pi_1(P_1) \cong
\pi_1(P_2)$.
\end{example}

\begin{example}(\textit{Symplectic fibrations})\label{ex:fibr}

It follows from the previous example that every simply-connected
symplectic manifolds is Morita equivalent to a point. Similarly,
if $(Q,\Pi)$ is a Poisson manifold with $\Pi=0$, then $Q$ is
Morita equivalent to any product $Q \times S$ where $S$ is a
simply-connected symplectic manifold. In fact, $Q\times S
\stackrel{\pr_1}{\leftarrow} T^*Q\times S
\stackrel{\pr_2}{\rightarrow} Q$ is a Morita bimodule, where
$\pr_1$ and $\pr_2$ are the natural projections.

More generally, let us assume that $P$ is
 a Poisson manifold whose symplectic foliation
is a fibration  $P\to Q$ with simply-connected fibres.
In general, there are obstructions to $P$ being Morita
equivalent to $Q$ \cite{Xu91}: $P$ is
Morita equivalent to $Q$ if and only if there exists a closed
2-form on $P$ which restricts to the symplectic form on each
fibre.
We will have more to say about ``fibrating'' Poisson manifolds and their Morita
invariants in Section \ref{subsec:moritainv}.
\end{example}

\begin{example}(\textit{Lie-Poisson structures})\label{ex:morita-lie-poiss}

Let us consider $\mathfrak{g}_1^*$ and $\mathfrak{g}_2^*$, the
duals of the Lie algebras $\mathfrak{g}_1$ and $\mathfrak{g}_2$,
equipped with their Lie-Poisson structures. Then
$\mathfrak{g}_1^*$ and $\mathfrak{g}_2^*$ are Morita equivalent if
and and only they are isomorphic. Indeed, suppose that
$\mathfrak{g}_1^* \stackrel{J_1}{\leftarrow} S
\stackrel{J_2}{\rightarrow} \mathfrak{g}_2^*$ is a Morita
bimodule, and let $\X=J_2^{-1}(0)$. A dimension count shows that
there exists $\mu \in \mathfrak{g}_1^*$ such that
$\X=J_1^{-1}(\mu)$. Since $S$ is a biprincipal bibundle for the
symplectic groupoids $\grd(\mathfrak{g}_i^*)=T^*G_i$, $i=1,2$, it
follows that $\X$ is a $(G_1,G_2)$-Morita bibundle. Therefore
$G_1$ and $G_2$ are isomorphic, and so are $\mathfrak{g}_1^*$ and
$\mathfrak{g}_2^*$.

This example also follows from the Morita invariants discussed in
Section \ref{subsec:moritainv}.

\end{example}

\begin{example}(\textit{Topologically stable structures on surfaces})\label{ex:TSS}

Let $\Sigma$ be a compact, connected, oriented surface equipped
with a Poisson structure $\Pi$ which has at most linear
degeneracies and whose zero set consists of $n$ smooth, disjoint,
closed curves, for $n\geq 0$. These are called
\textbf{topologically stable structures} (TSS) \cite{Rad02}.

Any two modular vector fields for $\Pi$ \cite{We97} coincide at
points where $\Pi$ vanishes, so the curves in the zero set carry a
natural orientation. We denote the zero set of $\Pi$, regarded as
an oriented 1-manifold, by $Z(\Sigma,\Pi)$. Two TSS $(\Sigma,\Pi)$
and $(\Sigma,\Pi')$ are \textbf{topologically equivalent} if there
is an orientation-preserving diffeomorphism $\psi:\Sigma \to
\Sigma'$ such that $\psi(Z(\Sigma,\Pi))=Z(\Sigma',\Pi')$. We
denote the equivalence class of $Z(\Sigma,\Pi)$ by
$[Z(\Sigma,\Pi)]$. This class can be represented by an oriented
labeled graph $\mathfrak{G}(\Sigma,\Pi)$: each vertex corresponds
to a 2-dimensional leaf of the structure, two vertices being
connected by an edge for each boundary zero curve they share; each
edge is oriented to point toward the vertex for which $\Pi$ is
positive with respect to the orientation of $\Sigma$. We then
label each vertex by the genus of the corresponding leaf.

 It turns out that
the topology of the zero set plus the modular periods (periods of
a modular vector field around the zero curves) completely
determine the Morita equivalence class of TSS \cite{BuRad02,BuWe}.
In fact, let us define a more elaborate graph
$\mathfrak{G}_T(\Sigma,\Pi)$, obtained from
$\mathfrak{G}(\Sigma,\Pi)$ by labeling each of its edges by the
modular period around the corresponding zero curve. Then two TSS
$(\Sigma,\Pi)$ and $(\Sigma',\Pi')$ are Morita equivalent if and
only if there is an isomorphism of labeled graphs
$\mathfrak{G}_T(\Sigma,\Pi)\cong \mathfrak{G}(\Sigma',\Pi')$. (It
follows from the results in \cite{CrFe01} that TSS are always
integrable.)

The classification of TSS up to Morita equivalence was preceded by
(and builds on) their classification up to orientation-preserving
Poisson diffeomorphisms by Radko \cite{Rad02}, who shows that the
topological class of the zero set and the modular periods,
together with a certain volume invariant (generalizing the
Liouville volume when the TSS is symplectic), form a complete set
of invariants.

\end{example}

\begin{remark}(\textit{Morita equivalence of presymplectic groupoids and
``momentum map theories''})

As we noted in Examples \ref{ex:hamilt} and
\ref{ex:cartdirmod}, hamiltonian spaces can be seen as modules
over Lie-Poisson structures on duals of Lie algebras, whereas
quasi-hamiltonian (or hamiltonian quasi-Poisson) manifolds are
modules over Cartan-Dirac structures on Lie groups. Thus, the
category of modules over an arbitrary (integrable) Poisson or
Dirac manifold can be regarded as the category of ``hamiltonian
spaces'' for some generalized ``momentum map theory''. Since Morita
equivalence establishes an equivalence of categories of modules,
it provides a precise notion of equivalence for
``momentum map theories'' and automatically implies the existence of
other invariants (such as classical intertwiner spaces--see Remark
\ref{rem:intertwiners}).

An extended notion of Morita equivalence for $\phi$-twisted presymplectic
groupoids (or, infinitesimally, $\phi$-twisted Dirac structures)
was developed by Xu in \cite{Xu03}. In Xu's  work, it is shown
that various known correspondences of ``momentum map theories''
can be described by appropriate Morita equivalences.
Examples include the equivalence between ordinary momentum maps
and momentum maps for actions of Poisson-Lie groups (taking values
in the dual group) \cite{Alek,GinzWe} and the one between
quasi-hamiltonian spaces for groups
and ordinary hamiltonian spaces for their loop groups \cite{AMM}.
An interesting feature of Morita equivalence for presymplectic groupoids is that
the bimodules are not simply a pair of modules structures which commute.

Besides relating ``momentum map theories'', Morita equivalence of
groupoids
plays a central role in certain approaches to geometric
quantization of these generalized hamiltonian spaces, where the usual
line
bundles are replaced by gerbes \cite{BeXuZh,LaXu}.
\end{remark}

\subsection{Picard groups}\label{subsec:picard}

Just as for algebras, there are Picard groupoids associated with
the categories $\catpoiss$ and $\SGrd$. In particular,
the isomorphism classes of Morita self-equivalences of a Poisson
manifold $P$  (resp. symplectic groupoid $\grd$) form a group
$\Pic(P)$ (resp. $\Pic(\grd)$), called the \textbf{Picard group}.
It follows from the discussion in the
previous section that $\Pic(P) = \Pic(\grd(P))$.

We now discuss some examples of ``geometric'' Picard groups; see
\cite{BuWe} for details.

\begin{example}(\textit{Picard groups of groups})\label{ex:picgrgr}

As we saw in Section \ref{subsec:moritagr}, geometric Morita
equivalences between groups are closely related to group
isomorphisms.  A closer analysis shows that the Picard group of a
group $G$ is naturally isomorphic to its group
$\Outaut(G):=\Aut(G)/\Inaut(G)$ of outer automorphisms.

It follows from Example \ref{ex:transitive} and the invariance of
Picard groups under Morita equivalence that, if $\grd$ is a
transitive  groupoid over $P$,
then $\Pic(\grd) \cong \Outaut(\grd_x)$, where $\grd_x$ is the isotropy group
at a point $x\in P$.  This isomorphism is natural, so the outer
automorphism groups attached to different points are all naturally
isomorphic to one another.
\end{example}

\begin{example}(\textit{Picard groups of symplectic
    manifolds})\label{ex:picsym}

Since, according to Example \ref{ex:sym}, the fundamental groupoid
of a connected symplectic manifold $P$ is Morita equivalent to any
of its fundamental groups $\pi_1(P,x)$, it follows from Example
\ref{ex:picgrgr} that, for such a manifold, $\Pic(P)$ is naturally
isomorphic to $\Outaut(\pi_1(P,x))$ for any $x$ in $P$.
\end{example}

The Picard group of a Poisson
manifold or symplectic manifold is also related to
 a group of outer automorphisms of
the manifold itself.
For a Poisson manifold $P$, let $\Aut(P)$ denote its group
of Poisson diffeomorphisms. There is a natural map
\begin{equation}\label{ex:geominc}
\inc:\Aut(P)\to \Pic(P),
\end{equation}
analogous to \eqref{eq:jmap}, which assigns to each $\psi\in
\Aut(P)$ the isomorphism class of the bimodule $P
\stackrel{t}{\longleftarrow}\grd(P) \stackrel{\psi^{-1}\circ
s}{\longrightarrow}{P}$. Any lagrangian bisection of $\grd(P)$
(which is the analogue of a group element) naturally
induces a Poisson diffeomorphism of $P$ that we call an
\textbf{inner automorphism}. It turns out that
$\ker(\inc)=\Inaut(P)$, the group of inner automorphisms of $P$,
just as in the algebraic setting discussed in Section
\ref{subsec:rings}.

The situation for symplectic groupoids is completely analogous \cite{BuWe}.

\begin{exercise}
Let $P$ be the $2n$ dimensional torus $\reals^{2n}/(2
\pi\integers)^{2n}$
with a symplectic
structure of the form $\frac{1}{2}\omega_{ij} d\theta^i \wedge d
\theta^j$, where $\omega$ is a nondegenerate antisymmetric matrix of
real constants.  Show that the Picard group of $P$ is independent of
the choice of $\omega$, while the subgroup of $\Pic(P)$ arising from
outer automorphisms of $(P,\omega)$ does depend on $\omega$.
\end{exercise}

\begin{exercise}
Compare $\Outaut(P)$ with $\Pic(P)$ when $P$ is the disjoint
union of several $2$-dimensional spheres, possible with different
symplectic areas.  Hint: use the theorems of Smale
\cite{Sm} and Moser \cite{mo:volume} to show that
every symplectomorphism of $S^2$ is inner.
\end{exercise}

There are  geometric versions of the maps \eqref{eq:projX} and
\eqref{eq:proj}.
Let  $\PSP$ be a Morita equivalence. If
$\mathcal{O}\subseteq P_2$ is a symplectic leaf, then
$J_1(J_2^{-1}(\mathcal{O}))$ is a symplectic leaf of $P_1$, and this
is a bijective correspondence between symplectic
leaves. So, for a Poisson manifold $P$, we have a map
\begin{equation}\label{eq:poisstatic}
\Pic(P) \to \Aut(\mathrm{Leaf}(P)),
\end{equation}
where $\mathrm{Leaf}(P)$ is the leaf space of $P$, analogous to
the map \eqref{eq:proj}. We define the \textbf{static Picard
group} $\SPic(P)$ of $P$ as the kernel of \eqref{eq:poisstatic},
i.e., the self-Morita equivalences inducing the identity map on
the leaf space.   Note that functions on the leaf space constitute
the center of the Poisson algebra of functions on $P$, hence the
terminology analogous to that for algebras.

\begin{example}(\textit{Zero Poisson structures})

As we saw in Example \ref{ex:zero}, in this case
$\grd(P)=T^*P$, and $\Pic(P)=\Pic(T^*P)$. Since
$\mathrm{Leaf}(P)=P$, \eqref{eq:poisstatic} implies that each
self-Morita bimodule $S$ induces a diffeomorphism
$\psi$ of $P$. So composing $S$ with $\psi^{-1}$ defines an
element of the static Picard group $\SPic(P)$.
A direct computation shows that the map \eqref{eq:poisstatic} is split by
the map $\Aut(P) \to \Pic(P)$  \eqref{ex:geominc}, hence
$$
\Pic(P) = \Diff(P) \ltimes \SPic(P),
$$
in complete analogy with \eqref{eq:piccomm}. Bimodules in $\SPic(P)$ are of the form
\begin{equation}\label{eq:Spp}
\xymatrix { {S} \ar@<-0.5ex>[d]_{p}\ar@<0.5ex>[d]^{p}\\
P}
\end{equation}
so each fibre $p^{-1}(x)$ is lagrangian and simply-connected;
moreover, the fact that $p$ is a \textit{complete} Poisson map
implies that the  $p$-fibres are complete with respect to their
natural affine structure.

Since $P\leftarrow S \to P$ is a Morita bimodule, the $p$-fibres
are isomorphic to the fibres of the symplectic groupoid target map
$T^*P\to P$, so they are contractible. As a result, there exists a
cross section $P \to S$, which implies that there is a
diffeomorphism $S\cong T^*P$ preserving the fibres
\cite[Sec.~3]{BuWe}. Hence, in order to characterize a bimodule
\eqref{eq:Spp}, the only remaining freedom is on the choice of
symplectic structure on $T^*P$. It turns out that the most general
possible symplectic structure on $T^*P$ with respect to which the
fibres of $T^*P\to P$ are lagrangian and complete is of the form:
$$
\omega + p^*B,
$$
where $\omega$ is the canonical symplectic form on $T^*P$ and $B$
is a closed 2-form on $P$ (a ``magnetic'' term). One can show that
two such bimodules are isomorphic if and only if $B$ is exact.
Hence
\begin{equation}
\SPic(P) \cong H^2(P,\mathbb{R}),
\end{equation}
and
\begin{equation}\label{ex:piczero}
\Pic(P) \cong \Diff(P) \ltimes H^2(P,\mathbb{R}),
\end{equation}
where the semi-direct product is with respect the natural action of $\Diff(P)$
on $H^2(P,\mathbb{R})$ by pull back. The reader can find the details in
\cite[Sec.~6.2]{BuWe}.

\begin{remark} (\textit{An intriguing resemblance})\label{rem:intriguing}

Recall from Example \ref{ex:picfn} that the Picard group of the
algebra $C^\infty(P)$ (which can be seen as a trivial quantization
of $(P,\Pi)$, if $\Pi=0$) is $\Diff(P)\ltimes H^2(P,\mathbb{Z})$.
Is there a theorem relating classical and quantum Picard groups
which would explain the similarity between this fact and
\eqref{ex:piczero}?
\end{remark}

\end{example}


\subsection{Fibrating Poisson manifolds and Morita invariants}
\label{subsec:moritainv}

In this section, we will discuss ``rigidity'' aspects of geometric
Morita equivalence.  As we saw in Theorem \ref{thm:equivcatpoiss},
Morita equivalence preserves categories of ``geometric
representations''.  We point out a few other invariants, some of which
have already appeared in previous sections.
\begin{enumerate}

\item As shown in Example \ref{ex:sym}, the Morita equivalence class
of a symplectic manifold is completely determined by the isomorphism
class of its fundamental group;

\item As remarked in Section \ref{subsec:picard}, Morita
equivalence induces a one-to-one correspondence of symplectic
leaves, which is a diffeomorphism whenever the leaf spaces are
smooth; moreover, corresponding symplectic leaves are themselves
Morita equivalent \cite{BuRad02,CrFe02} and have isomorphic
transverse Poisson structures \cite{We83};

\item Morita equivalence preserves first Poisson cohomology groups
\cite{Cra,GinzLu92}, and modular classes \cite{Cra,We97};

\item The monodromy groups and isotropy Lie algebras are Morita
invariant \cite{CrFe02}.

\end{enumerate}

As remarked in \cite{CrFe02}, all the invariants listed above turn out
to be preserved by a notion of equivalence which is much weaker than
Morita equivalence, called \textbf{weak Morita equivalence}, which
does not require the integrability of Poisson manifolds.  We do not
 know any Morita invariant which is not a weak Morita invariant.

By 1. above, the only Morita invariant of a connected symplectic
manifold is its fundamental group.  For a disjoint union of
symplectic components, it is the unordered list of fundamental
groups which counts; in particular, if all the components are
simply connected, the number of components is a complete
invariant.  In this section, we will see that the Morita invariant
structure is much richer for a Poisson manifold which is a smooth
family of (diffeomorphic) symplectic manifolds.

We will say that a Poisson manifold $P$ is \textbf{fibrating}
if its symplectic leaves are the fibres of a smooth
locally trivial fibration from $P$ to $\mathrm{Leaf}(P)$.
Here, locally triviality is meant in the differentiable rather than
symplectic sense; in fact, it is the variation in symplectic structure
from fibre to fibre which will concern us.

When $P$ is fibrating, the fibrewise homology groups
$H_2(\fib,\mathbb{Z})$ form a locally trivial bundle of abelian groups
over $\mathrm{Leaf}(P)$.  Pairing with the
fibrewise symplectic structure gives a map
$H_2(\fib,\mathbb{Z}) \to \mathbb{R},
$
which encodes the variation of the symplectic cohomology class from fibre to
fibre.   The derivative of this map with respect to the base point in
$\mathrm{Leaf}(P)$ gives rise to a map
$
\nu: H_2(\fib,\mathbb{Z}) \to T^*\mathrm{Leaf}(P).
$

The map $\nu$ vanishes on torsion elements of $H_2(\fib,\mathbb{Z})$,
so its image is a family of embedded abelian groups
 in the fibres of $T^*\mathrm{Leaf}(P)$, called the
\textbf{variation lattice}  of $P$.   Dazord \cite{Da} proves that,
if $P$ is integrable and has simply connected fibres,
the variation lattice must be topologically closed with constant rank,
having local bases of closed 1-forms.  Failure of the variation
lattice to have these properties provides an obstruction to
integrability which was
extended to general Poisson manifolds in \cite{CrFe02}.

A nice application of the variation lattice is to the study of the Picard
groups of the duals of Lie algebras of compact groups \cite{BuWe}, in
which the lattice imparts a flat affine structure to the regular part
of the symplectic leaf space.

\begin{example}\label{ex:noninteg2}(\textit{Nonintegrable Poisson structures
revisited})

Let us again consider $P = \mathbb{R} \times S^2$ from Example
\ref{ex:noninteg}, with Poisson
structure $(1/f(t))\Pi_{S^2} \times 0$, $f(t)>0$. The area
of the symplectic leaf over $t \in
\mathbb{R}$ is $4\pi f(t)$. The variation lattice is spanned by
$4\pi f'(t)dt$, so it has constant rank if and only if $f'(t) \equiv 0$
or $f'(t)$ is not zero for all $t$.
\end{example}

If $P_1$ and $P_2$ are Morita equivalent fibrating Poisson manifolds
with simply connected leaves, then the induced diffeomorphism $
\mathrm{Leaf}(P_1) \to \mathrm{Leaf}(P_2) $ preserves the variation
lattice; this can be seen as a special case of 4. above.  So, although
Morita equivalence does not determine the fibrewise symplectic
structures, it is sensitive to how symplectic structures vary from
fibre to fibre.  This sensitivity leads to the following example
\cite{Xu91}, of Poisson manifolds which are representation equivalent
but not Morita equivalent (see Remark  \ref{rem:categ}).

\begin{example}(\textit{Representation equivalence
vs.~Morita equivalence})\label{rem:repvsmorita}

Consider $(0,1) \times S^2$ with Poisson
structures $\Pi_1$ and $\Pi_2$ determined by the fibrewise
symplectic structures $(1/t)\Pi_{S^2}$ and $(1/2t)\Pi_{S^2}$,
respectively. Their variation lattices are spanned by $4\pi dt$
and $8 \pi dt$, respectively. Since there is no diffeomorphism of
$(0,1)$ taking $dt$ to $2 dt$, these structures cannot be Morita
equivalent. Note however that these structures are representation
equivalent: representations of $\Pi_1$ and $\Pi_2$ can be
interchanged by dividing or multiplying the symplectic form on the
realizations by 2.
\end{example}

\begin{remark} (\textit{A complete invariant?})

Xu \cite{Xu91} shows that the leaf space with its variation lattice
completely determines the Morita equivalence class of a fibrating
Poisson manifold for which the symplectic leaves are simply connected
and form a differentiably {\em globally} trivial fibration.  It does
not seem to be
known whether this result persists without the global triviality
assumption.  The attempt to attack this problem by ``gluing'' together
applications of the known case to local trivializations seems to lead
to the problem of computing the static Picard group of a fibrating Poisson
manifold.
\end{remark}

To extend the discussion above to the case where the leaves are not
simply connected, it seems that the variation lattice should be
replaced by its ``spherical'' part, obtained by replacing $H_2(\fib,\mathbb{Z})$
by the subgroup consisting of the spherical classes, i.e. the image of
the Hurewicz homomorphism from the bundle  $\pi_2(\fib)$ of homotopy groups.  This spherical
variation lattice is very closely related to the monodromy groups in
\cite{CrFe02}.  Details in this case should be interesting to work out,
particularly when the symplectic leaf fibration
is not globally trivial.

\begin{remark} (\textit{Noncompact fibres})

If the leaves of a fibrating Poisson manifold are compact, Moser's
theorem \cite{mo:volume} implies that the variation lattice actually
measures how the isomorphism class of the symplectic structure varies
from leaf to leaf.  If the leaves are noncompact, e.g. if they are
discs in $\reals^2$, then their area can vary without this being
detected by any Morita invariant.  Is there another notion of Morita
invariance which would detect the variation from fibre to fibre of
symplectic volume or other invariants, such as capacities?
\end{remark}

\begin{remark}  (\textit{Morita equivalence for nonintegrable Poisson manifolds})\label{rem:smorita}.

For a fibrating Poisson manifold which is nonintegrable, the variation
lattice still exists, so one might hope that it is still a Morita
invariant when the leaves are simply connected.  But there is no
Morita equivalence between such a manifold and itself, much less
another one.  To remedy this problem, we should extend the notion of
Morita equivalence to admit as bimodules smooth stacks which are not
manifolds, as we did for self-equivalences in Section \ref{subsec:sg}.
If we do this, then the variation lattice is indeed Morita invariant.
In particular, this shows that integrability is an invariant property
under this broadened notion of  Morita invariance.  Moreover, it turns
out that any ``S"-Morita equivalence between integrable Poisson manifolds
is given by a manifold, so that the integrable part of the Picard groupoid
remains unchanged.
It would be interesting to see how S-Morita equivalence is related to weak
Morita equivalence.
\end{remark}

\subsection{Gauge equivalence of Poisson structures}\label{subsec:gauge}

Let $P$ be a manifold, and let $\phi\in \Omega^3(P)$ be closed.
There is a natural way in which \textit{closed} 2-forms on $P$ act on
$\phi$-twisted Dirac structures: if $B\in \Omega^2(P)$ is closed
and $L$ is a $\phi$-twisted Dirac structure on $P$, then we set
$$
\tau_B(L):= \{(X,\alpha+\widetilde{B}(X))\,|\, (X,\alpha) \in L\},
$$
which is again a $\phi$-twisted Dirac structure. We call this
operation on Dirac structures a \textbf{gauge transformation}
associated with a 2-form \cite{SeWe01}. (More generally, for an
arbitrary $B$, $\tau_B(L)$ is a $(\phi-dB)$-twisted Dirac
structure.) Geometrically, a gauge transformation changes a Dirac
structure $L$ by adding the pull-back of a closed 2-form to its
leafwise presymplectic form.

\begin{remark}(\textit{Gauge transformations and $B$-fields})

In a complete similar way, complex closed 2-forms act on complex
Dirac structures. If $B \in \Omega^2(P)$ is a \textit{real}
2-form, and $L$ is a generalized complex structure on $P$ (see
Remark \ref{rem:complex}), then one can show that $\tau_B(L)$ is
again a generalized complex structure, and this operation is
called a  \textbf{$B$-field transform} \cite{Gualt,Hitchin}.

\end{remark}

If $\Pi$ is a Poisson structure on $P$, then changing it by a
gauge transformation will generally result in a Dirac structure
which is no longer Poisson. In fact, if $B\in \Omega^2(P)$ is
closed, then $\tau_B(L_\Pi)$ is a Poisson structure if and only if
the bundle map
\begin{equation}
\id + \widetilde{B}\widetilde{\Pi}: T^*P \to T^*P
\end{equation}
is invertible. In this case, the resulting Poisson structure is
the one associated with the bundle map
$$
\widetilde{\Pi}(\id + \widetilde{B}\widetilde{\Pi})^{-1}:T^*P \to TP,
$$
and we denote it by $\tau_B(\Pi)$.

Let $(P,\Pi)$ be a fibrating Poisson manifold, as in Section
\ref{subsec:moritainv}. Since a gauge transformation adds the
pull-back of a closed 2-form on $P$ to the symplectic form on each
fibre, the cohomology classes of fibrewise symplectic forms may
change in this operation; however, the way they vary from fibre to
fibre does not. This suggests that gauge transformations preserve
the Morita equivalence class of $(P,\Pi)$. In fact, this holds in
complete generality \cite{BuRad02}:
\begin{theorem}\label{thm:gaugemorita}
Gauge equivalence of integrable Poisson structures implies Morita
equivalence.
\end{theorem}

Since gauge transformations do not change the foliation of a
Poisson structure, there is no hope that the converse of Theorem
\ref{thm:gaugemorita} holds, since even Poisson diffeomorphic
structures may have different foliations. We call two Poisson
manifolds $(P_1,\Pi_1)$ and $(P_2,\Pi_2)$ \textbf{gauge equivalent
up to Poisson diffeomorphism} if there exists a Poisson
diffeomorphism $\psi:(P_1,\Pi_1)\to (P_2,\tau_B(\Pi_2))$ for some
closed 2-form $B\in \Omega^2(P_2)$. It clearly follows from
Theorem \ref{thm:gaugemorita} that if two integrable Poisson
manifolds are gauge equivalent up to a Poisson diffeomorphism,
then they are Morita equivalent. The following properties are
clear:

\begin{enumerate}
\item Two symplectic manifolds are gauge equivalent up to Poisson diffeomorphism if and only if
they are symplectomorphic;

\item If two Poisson manifolds are gauge equivalent up to Poisson
diffeomorphism, then they have isomorphic foliations (though generally
different leafwise symplectic structures);

\item The Lie algebroids associated with gauge equivalent Dirac structures are isomorphic \cite{SeWe01};
as a result, two Poisson manifolds which are gauge equivalent up to Poisson diffeomorphism
have isomorphic Poisson cohomology groups.
\end{enumerate}

A direct comparison between the properties above and the Morita
invariants listed in Section \ref{subsec:moritainv} suggests that
Morita equivalence should be a weaker notion of equivalence.
Indeed, two nonisomorphic symplectic manifolds with the same
fundamental group are Morita equivalent, but not gauge equivalent
up to Poisson diffeomorphism. In \cite[Ex.~5.2]{BuRad02},
one can also find examples of Morita equivalent Poisson structures
on the \textit{same} manifold which are not gauge equivalent up to
Poisson diffeomorphism by finding nonequivalent symplectic
fibrations with diffeomorphic total space and base (and using
Example \ref{ex:fibr}). Nevertheless, there are interesting
classes of Poisson structures for which both notions of
equivalence coincide, such as the topologically stable structures
of Example \ref{ex:TSS} \cite{BuRad02,BuWe}.

\begin{remark}(\textit{Gauge transformations and Morita equivalence of
quantum algebras})\label{rem:backalg}

As mentioned in Remark \ref{rem:formality}, gauge transformations
associated with \textit{integral} 2-forms define an action of
$H^2(P,\mathbb{Z})$ on formal Poisson structures on $P$ which can
be ``quantized'' (via Kontsevich's quantization \cite{Kont}) to
Morita equivalent deformation quantization algebras.

On the other hand, gauge transformations of translation-invariant
Poisson structures on tori are particular cases of the linear
fractional transformations \eqref{eq:fractional}, which quantize,
according to Theorem \ref{thm:MEqtori}, to strongly Morita
equivalent quantum tori.  As we already mentioned in Remark
\ref{rem:intriguing},
it would be very interesting to have a
unified picture relating Morita equivalence of quantum algebras to
geometric Morita equivalence.
\end{remark}

\section{Geometric representation equivalence}\label{sec:geomrep}

In Section \ref{sec:geometric},
we considered the category of $P$-modules (i.e. complete symplectic
realizations) over a Poisson manifold $P$, the geometric
analogue of the category of left modules over an algebra. We
observed in Remark \ref{rem:repvsmorita} that,
unlike the category of representations of an algebra, this category
does not determine the Morita equivalence class of $P$.
In this section, we will discuss refinements of the notion of
category of representations of a Poisson manifold in order to remedy
this defect.

\subsection{Symplectic torsors}

The first refinement we discuss is motivated by the theory of differentiable
stacks \cite{BeXu,Metzler,Pr}.

Given a Lie groupoid $\grd$, let $B\grd$ denote the category of
(left) $\grd$-torsors.
If two Lie groupoids $\grd_1$
and $\grd_2$ are Morita equivalent, then the natural functor $B\grd_1
\to B\grd_2$ induced by any Morita bibundle establishes an equivalence
of these categories.

However, to recover the Morita equivalence class of $\grd$ from $B\grd$,
one needs to consider another piece of information: the natural ``projection'' functor
$B\grd \to \mathcal{C}$, where $\mathcal{C}$ denotes the category of
smooth manifolds, which assigns to a $\grd$-torsor $S \to M$ the manifold $M=S/\grd$.
The category $B\grd$ \textit{together} with this projection functor is an example
of a differentiable stack. Taking this extra structure into account,
one defines $B\grd_1$ and $B\grd_2$ to be \textbf{isomorphic}
if there is an equivalence of categories $B\grd_1 \to B\grd_2$
commuting with the respective ``projections'' into $\mathcal{C}$.

It is clear that a functor induced by a Morita bibundle
establishes an isomorphism of stacks of torsors. It turns out that
the converse is also true:
if $B\grd_1$ and $B\grd_2$ are isomorphic in this refined sense,
then the Lie groupoids $\grd_1$ and $\grd_2$ are Morita equivalent.
As we will see, much of this discussion can
be adapted to the context of Poisson manifolds and symplectic groupoids.

Let $P$ be an integrable  Poisson manifold.
A  \textbf{symplectic $P$-torsor} is a complete
symplectic
realization
 $J:S \to P$ with the additional property that the induced left action of the
symplectic groupoid $\grd(P)$ on $S$ is principal. Note that, in this case,
the manifold $M=S/\grd(P)$ has a natural Poisson structure.  (As with
the regular bimodules in
Remark \ref{rem:regbim}, we can also describe
symplectic torsors without reference to groupoids.
$J:S \to P$ should be a surjective submersion, and the symplectic
orthogonal leaves to the $J$-fibres should be simply-connected and form a
simple foliation.)

Instead of considering the category of all complete symplectic realizations
over a Poisson manifold $P$, let us consider the category $BP$ of symplectic
$P$-torsors, as we did for
Lie groupoids.
If we restrict the morphisms in $BP$ to symplectomorphisms,
then there is a well-defined ``projection'' functor
$BP \to \mathcal{C}_{Pois}$, where $\mathcal{C}_{Pois}$ denotes the
category of Poisson manifolds, with \textit{ordinary} (invertible)
Poisson maps as
morphisms.
As in the case of Lie groupoids,
we refine the notion of isomorphism of categories to
include the ``projection'' functors: $BP_1$ and $BP_2$ are \textbf{isomorphic}
if there is an equivalence of categories $BP_1 \to BP_2$ commuting with
the projections $BP_i\to \mathcal{C}_{Pois}$, $i=1,2$.
In this setting, it is also clear that a Morita equivalence of $P_1$ and $P_2$
induces an isomorphism between $BP_1$ and $BP_2$. The following is a natural question.

\begin{quote}
If $BP_1$ and $BP_2$ are isomorphic, must $P_1$ and $P_2$ be
Morita equivalent Poisson manifolds?
\end{quote}

In Remark \ref{rem:repvsmorita}, we saw that the Poisson manifolds
$P_1=((0,1)\times S^2,(1/t)\Pi_{S^2})$ and $P_2=((0,1)\times S^2,(1/2t)\Pi_{S^2})$
are not Morita equivalent, but there is an equivalence of categories $BP_1 \to BP_2$.
However, this equivalence does \textit{not} commute with the ``projection'' functors,
so it is not an isomorphism in the refined sense.  Thus there is some
hope that the answer to the question above is ``yes,'' though we do
not yet have a complete proof.

\subsection{Symplectic categories}

The next approach to find a ``category of representations'' that
determines the Morita equivalence class of a Poisson manifold
is based on the
notion of ``symplectic category''. One can
think of it as the classical limit of the usual notion of abelian category, in the sense
that the vector
spaces (or modules) of morphisms in the theory of abelian categories
are replaced by symplectic manifolds.
Notice that we are
referring to \textit{a} ``symplectic category'', rather than \textit{the}
symplectic ``category'' of \cite{We82}.  From now on, we will drop the
quotation marks when referring to the new notion.

In a \textbf{symplectic category}, one has a class of objects, and,
for any two objects $A$ and $B$, a symplectic manifold, denoted by
$\mathrm{Hom}(A,B)$, which plays the role of the space of morphisms from $B$ to $A$.
Given three objects $A$, $B$ and $C$, the ``composition operation''
$\mathrm{Hom}(A,C) \leftarrow \mathrm{Hom}(A,B) \times \mathrm{Hom}(B,C)$
is  a lagrangian submanifold
$$
L_{ABC}\subset
\mathrm{Hom}(A,C) \times \overline{\mathrm{Hom}(A,B)} \times \overline{\mathrm{Hom}(B,C)}.
$$
This may \textit{not} be the graph of a map, but just a canonical relation, so we will
refer to it as the \textbf{composition relation}.
So, unlike in ordinary categories, $\mathrm{Hom}(A,B)$ should not be
thought of as a set of points.  Instead,
certain lagrangian submanifolds of $\mathrm{Hom}(A,B)$ will play the role
of ``invertible elements'', so that we can talk about ``isomorphic'' objects.
In other words, the guiding principle is to think of \textit{a}
symplectic category as a category in \textit{the} symplectic ``category.''

A \textbf{functor} between symplectic categories
should consist of a map $F$ between objects together with symplectic maps
$\mathrm{Hom}(A,B)\to \mathrm{Hom}(F(A),F(B))$, so that the induced map
from $$\mathrm{Hom}(A,C) \times {\mathrm{Hom}(A,B)} \times {\mathrm{Hom}(B,C)}$$
to $$\mathrm{Hom}(F(A),F(C)) \times {\mathrm{Hom}(F(A),F(B))}
\times {\mathrm{Hom}(F(B),F(C))}$$
maps the composition relation $L_{ABC}$ to $L_{F(A)F(B)F(C)}$.
It is also natural to require that
if $\mathrm{Hom}(A,B)$ contains ``invertible elements'',
then so does $\mathrm{Hom}(F(A),F(B))$.

If $\mathcal{S}$ and $\mathcal{S}'$ are symplectic categories,
then a functor $\mathcal{S} \to \mathcal{S}'$ is an
\textbf{equivalence of symplectic categories} if for any object $A'$ in $\mathcal{S}'$,
there exists an object $A$ such that $F(A)$ and $A$ are ``isomorphic''
(in the sense that $\mathrm{Hom}(F(A),A')$ contains an ``invertible
element''), and  the maps $\mathrm{Hom}(A,B)\to \mathrm{Hom}(F(A),F(B))$ are
symplectomorphisms.

We have not answered some questions about symplectic categories which
arise naturally.  Is $\mathrm{Hom}(A,A)$ always a symplectic
  groupoid?  If not, what are sufficient conditions?  Is there always
  a ``base'' functor from a given symplectic category to the category of Poisson
  manifolds and Morita morphisms?   Nevertheless, we can still discuss
  interesting examples, such as the one which follows.

\subsection{Symplectic categories of representations}

In the theory of abelian categories, a model example is the
category of modules over a ring (for instance, the group ring of a group,
in which case we have a category of representations). The morphisms are
module homomorphisms (or intertwining operators in the case of
representations).  The symplectic analogue of this example is
the ``symplectic category'' of representations of a Poisson manifold,
in which objects are symplectic realizations and
spaces of morphisms are the classical intertwiner spaces.

To avoid smoothness issues, we will be more restrictive and define
the \textbf{representation category} of a Poisson manifold
$P$ to be the symplectic category in which
the objects are symplectic $P$-torsors $S \to M$ which are
$(P,M)$-Morita equivalences,
and the morphism spaces are classical intertwiner spaces,
$
\mathrm{Hom}(S_1,S_2):= \overline{S_2}*_{P}{S_1}.
$

Composition relations are given by
\begin{equation}\label{eq:comprel}
L_{S_1S_2S_3}:= \{ ([(z,x)],[(y,x)],[(z,y)]) \}\subset
\overline{S_3}*_{P}{S_1} \times {S_2}*_{P}\overline{S_1}
\times{S_3}*_{P} \overline{S_2}
\end{equation}
where $[(a,b)] \in \overline{S'}*_P {S}$ denotes the image of
$(a,b)\in \overline{S'}\times_P S$
under the natural projection.

\begin{exercise}
Check that the composition relation  \eqref{eq:comprel}
is a lagrangian submanifold.  (Hint: first prove it when $P$ is a
point, then
use coisotropic reduction for the general case.)
\end{exercise}

Note that $\mathrm{Hom}(S,S)=\overline{S}*_{P}{S}$
is symplectomorphic to the
symplectic groupoid $\grd(M)$, where $M=S/\grd(P)$.

\begin{exercise}
Show that the composition relation in $\mathrm{Hom}(S,S)=\grd(M)$,
where $M=S/\grd(P)$, is just the graph of the groupoid
multiplication $ \grd(M)\leftarrow \grd(M)\times_{(s,t)} \grd(M) $ in
$\grd(M)\times \overline{\grd(M)} \times {\grd(M)}$.
\end{exercise}

Finally, we define ``invertible elements''
in $\mathrm{Hom}(S_1,S_2)=\overline{S_2}*_{P}{S_1}$
to be those lagrangian submanifolds which are the reductions
of graphs of isomorphisms
of symplectic realizations $S_1\to S_2$ via the coisotropic
submanifold $\overline{S_2}\times_{P}{S_1}$ of
$\overline{S_2}\times {S_1}$.
In particular, two symplectic
realizations
are ``isomorphic'' in the representation category of $P$ if and only if they
are isomorphic in the usual sense.

\begin{proposition}
Two Poisson manifolds are Morita equivalent if and only if they have
equivalent representation categories.
\end{proposition}

\begin{proof}

Suppose that $P_1$ and $P_2$ are Morita equivalent, and let $\X$ be a $(P_1,P_2)$-Morita
bimodule. Let $\mathcal{S}(P_i)$ denote the representation category of $P_i$, $i=1,2$.
Then $\X$ induces an equivalence of symplectic categories
$\mathcal{S}(P_2) \to \mathcal{S}(P_1)$:
at the level of objects, a Morita bimodule
$P_2 \leftarrow S \to M$ is mapped to the Morita bimodule
$P_1 \leftarrow \X*_{P_2}S \to M$; if $P_1\leftarrow S' \to M'$ is an object in
$\mathcal{S}(P_1)$, then ${\X}*_{P_1}\overline{S'}$ is an object in $\mathcal{S}(P_1)$
such that $S'$ and $\X*_{P_2}\overline{\X}*_{P_1}S'$ are isomorphic;
at the level of morphisms, because $\overline{S_2}*_{P_2}{S_1}\cong
\overline{S_2}*_{P_2}\overline{\X}*_{P_1}\X *_{P_2}{S_1}$, we have a natural
symplectomorphism
$$
\mathrm{Hom}(S_1,S_2)\cong
\mathrm{Hom}(\X*_{P_1}S_1,\X*_{P_1}S_2).
$$

Conversely, suppose that $F: \mathcal{S}(P_2)\to \mathcal{S}(P_1)$
is an equivalence of symplectic categories, and let $P_2\leftarrow S\to M$
be an object in $\mathcal{S}(P_2)$.
Then there is a symplectomorphism from $\mathrm{Hom}(S,S)=\grd(M)$ to
$\mathrm{Hom}(F(S),F(S))=\grd(M')$, where $M'=F(S)/\grd(P_1)$. Since this
symplectomorphism preserves the composition relation, it is a symplectic groupoid
isomorphism. In particular, $M$ and $M'$ are isomorphic as Poisson manifolds,
which implies that $F(S)$ is a $(P_1,M)$-Morita bimodule. If we take
$S=\grd(P_2)$, then $M=P_2$ and $F(S)$ is a $(P_1,P_2)$-Morita bimodule.
\end{proof}

The equivalence in Remark \ref{rem:repvsmorita} does not preserve
intertwiner spaces (their symplectic structures are related by a
factor of 2), so it does not contradict the result above.

\begin{footnotesize}

\end{footnotesize}

\end{document}